\documentclass[twoside]{irmaems}
\usepackage{amssymb} 
\usepackage{amsmath} 
\usepackage{latexsym}
\usepackage{graphicx}
\usepackage{mathrsfs}

\renewcommand\theequation{\arabic{section}.\arabic{equation}}

\theoremstyle{definition} 

 \newtheorem{definition}{Definition}[section]
 \newtheorem{remark}[definition]{Remark}

\theoremstyle{plain}      

 \newtheorem{theorem}[definition]{Theorem}
 \newtheorem{corollary}[definition]{Corollary}
 
 \newtheorem*{example}{Example}

\newtheorem*{conjecture}{Conjecture}

\renewcommand{\vec}[1]{\mathbf{#1}}

\begin{document}

\title{Mirzakhani's recursion formula on Weil-Petersson volume and applications}

\author{Yi Huang
}

\address{
Department of mathematics and statistics,
University of Melbourne\\
University of Melbourne, Victoria 3010, Australia\\
email: {\tt huay@ms.unimelb.edu.au}
}

\maketitle

\begin{abstract}
We give an overview of the proof for Mirzakhani's volume recursion for the Weil-Petersson volumes of the moduli spaces of genus $g$ hyperbolic surfaces with $n$ labeled geodesic boundary components, and her application of this recursion to Witten's conjecture and the study of simple geodesic length spectrum growth rates.
\end{abstract}



\tableofcontents   

\section{Introduction}\label{backgroundsec}

Let $S_{g,n}$ denote a (connected) surface with genus $g$ and $n$ punctures labeled $1$ to $n$, so that the Euler characteristic $\chi(S_{g,n})=2-2g-n$ is negative. The moduli space $\mathcal{M}(S_{g,n})$ is a real $(6g-6+2n)$-dimensional orbifold, whose points represent isometry classes of complete finite-area hyperbolic metrics on $S_{g,n}$.

The moduli space $\mathcal{M}(S_{g,n})$ has a natural symplectic structure given by the Weil-Petersson symplectic form $\omega_{S_{g,n}}$, and the volume form obtained by taking the top exterior product
\begin{align*}
\Omega_{S_{g,n}}:=\tfrac{1}{(3g-3+n)!}\underbrace{\omega_{S_{g,n}}\wedge\ldots\wedge\omega_{S_{g,n}}}_{3g-3+n\text{ terms}}
\end{align*}
has finite volume. In particular, for the once-punctured torus $S_{1,1}$, Wolpert showed that this volume is $\frac{1}{6}{\pi^2}$ in two different ways \cite{wolpert1, wolpert2}:
\begin{enumerate}
\item
explicitly computing the volume of a fundamental domain for $\mathcal{M}(S_{1,1})$;
\item
identifying the Weil-Petersson form as a cohomology class, thus relating intersection numbers on $\mathcal{M}(S_{1,1})$ with its Weil-Petersson volume.
\end{enumerate}

The four-punctured sphere $S_{0,4}$ case was similarly derived from this first volume computation. However, it was not until Penner's work in \cite{penner} that a further Weil-Petersson volume was obtained. Penner gave a fatgraph-based\footnote{Also known as ribbon-graphs to those more familiar with Kontsevich's work.} cell-decomposition of $\mathcal{M}(S_{g,n})$ and used this to explicitly obtain that the volume of $\mathcal{M}(S_{1,2})$ is $\frac{1}{4}\pi^4$ \cite[Thm.~5.2.1]{penner}. Indeed, Penner described a general strategy for computing the volumes of the moduli space $\mathcal{M}(S_{g,n})$ for any punctured surface $S_{g,n}$. In practice, however, this is intractable because: 
\begin{enumerate}
\item
the number of cells in $\mathcal{M}(S_{g,n})$ grows quite quickly (by \cite[Thm.~B]{penner}, even for $n=1$ the growth rate tends to $\frac{(2g)! 6^{2g}}{(6g-3)\cdot e^{2g}}$) and
\item
the integral for the volume of each top-dimensional cell becomes difficult to exactly evaluate.
\end{enumerate}

Zograf, however, explicitly expressed the Poincar\'e dual of the Weil-Petersson form in terms of certain divisors on the Deligne-Mumford compactification locus of $\mathcal{M}(S_{0,n})$, and exploited the intersection number interpretation of the Weil-Petersson volume to obtain the following recursion formula for the volume $V_{0,n}$ of the moduli space of $n$-punctured spheres \cite{zograf1}:
\begin{gather*}
V_{0,n}=\frac{(2\pi^2)^{n-3}}{(n-3)!}\cdot v_{0,n}\text{, where }v_{0,3}=1\text{ and }\\
v_{0,n}=\frac{1}{2}\sum_{k=1}^{n-3}\frac{k(n-k-2)}{n-1}{n-4\choose k-1}{n\choose k+1 }v_{0,k+2}\cdot v_{0,n-k}\text{ for }n\geq4.
\end{gather*}
He expanded on this work in \cite{zograf2}, to derive the following recursion formula for the volume $V_{1,n}$ of the moduli space of $n$-punctured tori\footnote{Some readers may notice that Zograf's volume $V_{1,1}=\frac{\pi^2}{12}$ is only half of $\frac{\pi^2}{6}$. This is due to the fact that all once-punctured tori have hyperelliptic involutions, this is discussed in greater detail in Remark~\ref{hyperelliptic}.}
\begin{gather*}
V_{1,n}=\frac{(2\pi^2)^{n}}{n!}\cdot v_{1,n}\text{, where}\\
v_{1,n}=\frac{n}{24}\cdot v_{0,n+2}+\sum_{k=1}^{n-1}(n-k){n-1\choose k }{n\choose k-1}v_{1,k}\cdot v_{0,n-k+2}\text{ for }n\geq1.
\end{gather*}

Then came N\"a\"at\"anen and Nakanishi's work \cite{nn1,nn2} on computing the volumes of moduli spaces of once-punctured tori and four-punctured spheres with geodesic boundaries of length $b_j$ (and cone angle singularities $\theta_j$). Nakanishi and N\"a\"at\"anen found that the answer is a rational polynomial in $\pi^2$ and $b_j^2$ (or $\theta_j^2$): $\frac{1}{6}\pi^2+\frac{1}{24}b_1^2$ for the punctured torus and $2\pi^2 +\frac{1}{2}(b_1^2+b_2^2+b_3^2+b_4^2)$ for the four-punctured sphere. Their computation was also based on integrating over a fundamental domain, and thus was not easily generalizable.

This was the landscape prior to Mirzakhani's beautiful solution to the volume computation problem \cite{mirz1}. In Section~\ref{backgroundsec}, we give the necessary background on the Weil-Petersson geometry of moduli spaces; in Section~\ref{mcshanesection} we explain Mirzakhani's proof of her McShane identities --- a key element of her proof; in Section~\ref{volumesec} we give a schematic of Mirzakhani's proof strategy; in Section~\ref{volumesec} we use the McShane identity to give a sample computation of the volume of $\mathcal{M}(S_{1,2},(b_1,b_2))$; finally, in Section~\ref{applicationsec} we give applications of Mirzakhani's volume integration to proving Witten's conjecture and for specifying the growth rate of simple closed geodesics on a given hyperbolic surface.

\subsection{Preliminaries}

Unless otherwise specified, any surface $S$ that we consider is \emph{oriented}, \emph{hyperbolic}, \emph{finite area} and has either \emph{geodesic} or \emph{cusp} borders labeled from $1$ to $n$, where $n$ is the number of boundary components of $S$. In terms of defining the Teichm\"uller space, the moduli space and intermediate moduli spaces, we only require $S$ to be a topological surface. However, it is sometimes convenient to endow $S$ with a hyperbolic structure.

\subsubsection{Teichm\"uller space.}\index{Teichm\"uller space}
Let $S$ be a surface with genus $g$ and $n$ boundary components labeled $1$ to $n$, and let $\vec{b}=(b_1,\ldots,b_n)\in\mathbb{R}_{\geq0}^n$ be an $n$-tuple of positive real numbers, then the Teichm\"uller space $\mathcal{T}(S,\vec{b})$ is:
\begin{align*}
\left\{
\begin{array}{r|l}
 & X\text{ is a hyperbolic surface with labeled boundaries }\\
(X,f) & \text{of lengths $b_1,\ldots, b_n$, and }\\
& f:S\rightarrow X\text{ is a label-fixing homeomorphism}
\end{array}
\right\}/\sim_{\mathcal{T}},
\end{align*}
where $(X_1,f_1)\sim_{\mathcal{T}}(X_2,f_2)$ if and only if $f_2\circ f_1^{-1}:X_1\rightarrow X_2$ is isotopy equivalent to a isometry. We denote these equivalence classes by $[X,f]$ and refer to them as \emph{marked surfaces}\index{marked surface}.

We adopt the convention that length $b_i=0$ means that the $i$-th boundary component is a cusp, and we write $\mathcal{T}(S)$ for the Teichm\"uller space of a closed surfaces $S$. Teichm\"uller \cite{teich} showed that the Teichm\"uller space is homeomorphic to a $6g-6+2n$-dimensional open ball.

\subsubsection{Mapping class groups and moduli spaces.}

The group of boundary label-preserving homeomorphisms $\mathrm{Homeo}(S)$ of $S$ acts on the Teichm\"uller space by precomposition: given $h\in\mathrm{Homeo}(S)$,
\begin{align*}
h\cdot [X,f]:=[X,f\circ h^{-1}].
\end{align*}
Note that we precompose by the inverse of $h$ so that this is a left action.

Since marked surfaces are defined up to isotopy, the normal subgroup $\mathrm{Homeo}_0(S)$, consisting of all homeomorphisms isotopy equivalent to the identity map $\mathrm{id}:S\rightarrow S$, acts trivially on $\mathcal{T}(S,\vec{b})$. Thus, we define the \emph{mapping class group}\index{mapping class group}
\begin{align*}
\mathrm{Mod}(S):=\mathrm{Homeo}(S)/\mathrm{Homeo}_0(S).
\end{align*}
We refer to elements $[h]$ of $\mathrm{Mod}(S)$ as \emph{mapping classes}.

Mapping class groups are finitely presentable \cite[Thm.~5.3]{primer}, and act discretely on the Teichm\"uller space. This action is almost free, in the sense that isotropy groups are finite, and the resulting quotient orbifold
\begin{align*}
\mathcal{M}(S,\vec{b}):=\mathcal{T}(S,\vec{b})/\mathrm{Mod}(S)
\end{align*}
is referred to as the \emph{moduli space}\index{moduli space} of $S$. Note that each point of $\mathcal{M}(S,\vec{b})$ represents a distinct isometry class of hyperbolic surfaces homeomorphic to $S$ with boundary lengths equal to $\vec{b}$. Thus, as a set, we may identify $\mathcal{M}(S,\vec{b})$ with the set
\begin{align*}
\left\{
X\;
\begin{array}{|l}
X\text{ is a hyperbolic surface with}\\
\text{labeled boundaries of lengths }b_1,\ldots, b_n
\end{array}
\right\}/\sim_\mathcal{M},
\end{align*}
where $X_1\sim_{\mathcal{M}}X_2$ if and only if they are isometric. We denote the isometry class of $X$ by $[X]\in\mathcal{M}(S,\vec{b})$.

\subsubsection{Pairs of pants.}\index{pair of pants}

A \emph{pair of pants} is a hyperbolic surface with genus $0$ and $3$ geodesic boundary components (allowing for cuspidal boundaries).

\begin{theorem}[{\cite[Thm.~3.1.7]{buser}}]\label{theorem_pants}
For any 3-tuple $(b_1.b_2,b_3)\in\mathbb{R}_{\geq0}^3$, there is a unique hyperbolic pair of pants $S$ with labeled boundaries $\beta_1,\beta_2,\beta_3$ respectively of lengths $b_1,b_2$ and $b_3$.
\end{theorem}
 
Since every label-preserving self-homeomorphism on a pair of pants is homotopy equivalent to the identity map, we have:

\begin{corollary}
The moduli space $\mathcal{M}(S,(b_1,b_2,b_3))$ for a hyperbolic pairs of pants $S$ is precisely $\{[S]\}$.
\end{corollary}

The only simple closed geodesics on a pair of pants $S$ are the three geodesics $\beta_1,\beta_2,\beta_3$ constituting the boundary $\partial S$ of $S$. Moreover, there is a unique simple geodesic arc $\sigma_{i,j}$ joining boundaries $\beta_i$ and $\beta_j$.

\begin{figure}[h]
\begin{center}
\includegraphics[scale=0.2]{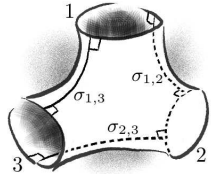}
\end{center}
\caption{A pair of pants with $\sigma_{1,2}$, $\sigma_{1,3}$ and $\sigma_{2,3}$}
\label{pg5}
\end{figure}

\begin{remark}
The three orthogeodesics\index{pair of pants!orthogeodesics} $\sigma_{1,2}=\sigma_{2,1}$, $\sigma_{1,3}=\sigma_{3,1}$ and $\sigma_{2,3}=\sigma_{3,2}$ cut $S$ up into two isometric right-angled hexagons. This tells us that there is an orientation-reversing isometry on $S$ which fixes these orthogeodesics pointwise and takes one hexagon to the other.
\end{remark}

\subsubsection{Fenchel-Nielsen coordinates.}\index{Fenchel-Nielsen coordinates}
Given a closed curve $\gamma$ on $S$, and a marked hyperbolic surface $[X,f]$, there is a unique closed geodesic homotopy equivalent to $f(\gamma)\subset X$. Denote this geodesic by $f_*\gamma$, then the curve $\gamma$ defines a positive valued function on Teichm\"uller space given by:
\begin{align*}
\ell_\gamma:\mathcal{T}(S,\vec{b})&\rightarrow\mathbb{R}_{> 0}\\
[X,f]&\mapsto\ell(f_*\gamma),
\end{align*}
where $\ell$ is the (geodesic) length function.

A \emph{pants decomposition}\index{pants decomposition} of a hyperbolic surface $S$, with genus $g$ and $n$ labeled boundaries, is a maximal collection $\{\gamma_1,\ldots,\gamma_{3g-3+n}\}$ of disjoint simple closed geodesics on $S$. As seen in Theorem~\ref{theorem_pants}, hyperbolic pairs of pants are uniquely determined by the lengths of their boundary geodesics. Thus, given a marked hyperbolic surface $[X,f]$, the pairs of pants obtained from cutting $X$ along $f_*\gamma_1,\ldots,f_*\gamma_{3g-3+n}$ are uniquely determined by the lengths $b_1,\ldots b_n$ and the length functions $\ell_i:=\ell_{\gamma_i}$ on Teichm\"uller space.

We can recover $X$ from a pants decomposition if we know how to glue its constituent pairs of pants. For each $\gamma_i$, the endpoints of the pants seams on the pairs of pants bordered by $\gamma_i$ allow us to keep track of this gluing with an element of $\mathbb{R}/\ell_i\mathbb{Z}$. Moreover, since any map $h: X\rightarrow X$ that fixes $f_*\gamma_i$ for $i=1,\ldots n$ is homotopy equivalent to Dehn twists of $X$ along the $\{f_*\gamma_i\}$, we may keep track of these Dehn twists (effectively lifting $\mathbb{R}/\ell_i\mathbb{Z}$ up to $\mathbb{R}$) and hence parameterize the entire Teichm\"uller space $\mathcal{T}(S,\vec{b})$. We denote the \emph{twist parameter} for $\gamma_i$ by $\tau_i$.

\begin{theorem}
Given a pants decomposition $\Gamma$ of a hyperbolic surface $S$ with genus $g$ and $n$ labeled boundaries of lengths $\vec{b}$, the Fenchel-Nielsen coordinates
\begin{align*}
\mathcal{FN}:\mathcal{T}(S,\vec{b})&\rightarrow \mathbb{R}_{>0}^{3g-3+n}\times\mathbb{R}^{3g-3+n}\\
[X,f]&\mapsto (\ell_1,\ldots,\ell_{3g-3+n},\tau_1,\ldots,\tau_{3g-3,n})
\end{align*}
is a real analytic homeomorphism.
\end{theorem}

\subsubsection{Weil-Petersson structure.}

Fix a pants decomposition $\{\gamma_1,\ldots,\gamma_{3g-3+n}\}$ on $S$, and consider the \emph{Weil-Petersson symplectic form}\index{Weil-Petersson!symplectic form} on $\mathcal{T}(S,\vec{b})$ given by:
\begin{align*}
\omega_{S,\vec{b}}:=\mathrm{d}\ell_1\wedge\mathrm{d}\tau_1+\ldots+\mathrm{d}\ell_{3g-3+n}\wedge\mathrm{d}\tau_{3g-3+n}.
\end{align*}
The Weil-Petersson form has a pants-decomposition independent formulation, and Wolpert showed that this $2$-form has the above expression \cite{wolpert3}. This means that $\omega_{S,\vec{b}}$ is invariant under the mapping class group action and descends to a symplectic 2-form on $\mathcal{M}(S,\vec{b})$. The Weil-Petersson volume obtained by integrating the top exterior product volume form\index{Weil-Petersson!volume form}
\begin{align*}
\Omega_{S,\vec{b}}:=\tfrac{1}{(3g-3+n)!}\underbrace{\omega_{S,L}\wedge\ldots\wedge\omega_{S,\vec{b}}}_{3g-3+n\text{ terms}}
\end{align*}
over $\mathcal{M}(S,\vec{b})$ is finite. Wolpert \cite{wolpert3} showed that $\omega_{S,\vec{b}}$ smoothly extends to a closed 2-form on the Deligne-Mumford compactification of $\mathcal{M}(S,\vec{b})$, and pairing its top exterior product with the fundamental class gives the Weil-Petersson volume. Indeed, the existence of Bers' constant\index{Bers' constant} $\mathit{Ber}(S)$ \cite[Thm.~5.1.2]{buser} implies that the finite volumed set
\begin{align*}
(0,\mathit{Ber}(S)]^{6g-6+2n}\subset  \mathcal{T}(S,\vec{b})
\end{align*}
contains a fundamental domain for the mapping class group. Although Bers' constant is not known exactly in general, there are upper bounds for $\mathit{Ber}(S)$ \cite[Thm.~4.8]{parlier} and so this does give us an explicit upper bound on the volume of $\mathcal{M}(S,\vec{b})$.

\subsection{Intermediate moduli spaces}\index{intermediate moduli space}
Since the Teichm\"uller space $\mathcal{T}(S,\vec{b})$ is contractible, it is the (orbifold) universal cover of $\mathcal{M}(S,\vec{b})$. Moreover, since the mapping class group is the group of deck transformations on $\mathcal{T}(S,\vec{b})$, subgroups of the mapping class group correspond to (connected) covering spaces of the moduli space. We informally refer to these covering spaces as \emph{intermediate moduli spaces}.

Now, given an ordered $m$-tuple $\Gamma=(\gamma_1,\ldots,\gamma_m)$ of simple closed geodesics on $S$, the mapping class group $\mathrm{Mod}(S)$ acts on $\Gamma$ diagonally to produce a collection 
\begin{align*}
\mathrm{Mod}(S)\cdot \Gamma=\left\{(h_*\gamma_1,\ldots,h_*\gamma_m)\mid h\in\mathrm{Mod}(S)\;\right\}
\end{align*}
of $m$-tuples of simple closed geodesics on $S$. The stabilizer 
\begin{align*}
\mathrm{Stab}(\Gamma):=\left\{h\in\mathrm{Mod}(S)\mid h_*\gamma_j=\gamma_j\text{ for }j=1,\ldots,m\;\right\}
\end{align*}
is a subgroup of the mapping class group, and corresponds to the following covering space:
\begin{align*}
\mathcal{M}(S,\Gamma,\vec{b}):=
\left\{
(X,f_*\alpha)\;
\begin{array}{|l}
[X,f]\text{ is a marked surface with}\\
\text{labeled boundaries of lengths }b_1,\ldots,b_n
\end{array}
\right\}/\sim_{\mathcal{M},\Gamma},
\end{align*}
where $(X_1,{f_1}_*\Gamma)\sim_{\mathcal{M},\Gamma}(X_2,{f_2}_*\Gamma)$ if and only if there is an isometry $h:X_1\rightarrow X_2$ such that $h_*({f_1}_*\gamma_j)={f_2}_*\gamma_j$ for $j=1,\ldots,m$. 

\subsubsection{Forgetful map.}
We refer to the covering map given by
\begin{align*}
\pi_{\Gamma}:\mathcal{M}(S,\Gamma,\vec{b})&\rightarrow\mathcal{M}(S,\vec{b})\\
[X,f_*\Gamma]&\mapsto[X]
\end{align*}
as the \emph{forgetful map}\index{forgetful map}, because $\pi_\Gamma$ forgets the geodesic $f_*\Gamma$ paired with $X$. By decorating surfaces $X$ with $m$-tuples of geodesics, we're able to define the following length function
\begin{align*}
\ell_\Gamma:\mathcal{M}(S,\Gamma,\vec{b})&\rightarrow\mathbb{R}_{>0}^m\\
[X,f]&\mapsto(\ell_{\gamma_1}([X,f]),\ldots,\ell_{\gamma_m}([X,f])).
\end{align*}

\subsubsection{Structure of intermediate moduli spaces.}
Let $\vec{c}=(c_1,\ldots,c_m)\in\mathbb{R}_{>0}^m$ denote the length $\ell_\Gamma([S,\Gamma])$ of $\Gamma$ on $S$, and let $S_1,\ldots,S_k$ denote the connected bordered hyperbolic surfaces resulting from cutting $S$ along $\Gamma$. Let us label/order the boundaries $S_1$ and denote the lengths of the borders which arise from cutting along $\Gamma$ by  
\begin{align*}
\vec{c}_{i}:=(c_{i,1},\ldots,c_{i,m_i})\in\{c_1,\ldots,c_m\}^{m_i}
\end{align*}
and the lengths of the other borders by
\begin{align*}
\vec{b}_i:=(b_{i,1},\ldots,b_{i,n_i})\in\{b_1,\ldots,b_n\}^{n_i}.
\end{align*}
Observe that $m_1+\ldots+m_k=2m$ and $n_1+\ldots+n_k=n$.

If $\Gamma$ consists of disjoint simple closed curves, then the multicurve length function $\ell_\Gamma$ maps surjectively onto its codomain $\mathbb{R}_{>0}^m$. The preimage $\ell_\Gamma^{-1}(\vec{c})$ of $\vec{c}$ consists of hyperbolic surfaces $X$ paired with an ordered $m$-tuple $f_*\Gamma$ of simple geodesics on $X$. Cutting $X$ along $f_*\Gamma$ results in hyperbolic surfaces $X_1,\ldots,X_k$ respectively homeomorphic to $S_1,\ldots, S_k$. To recover $[X,f]$ from the $\{X_i\}$, we only need to specify how these subsurfaces are glued together. Since the $X_i$ may vary over $\mathcal{M}(S_i,(\vec{c}_i,\vec{b}_i))$ and the gluing for $\gamma_j$ varies over $\mathbb{R}/c_j\mathbb{Z}$, we obtain that:
\begin{align}\label{levelset}
\ell_\Gamma^{-1}(\vec{c})
=(\mathbb{R}/c_{1}\mathbb{Z})\times\ldots\times(\mathbb{R}/c_{m}\mathbb{Z})
\times\prod_{i=1}^k\mathcal{M}(S_i,(\vec{c}_i,\vec{b}_i)) .
\end{align}
Note that the above identification~\eqref{levelset} holds for any $\vec{c}\in\mathbb{R}_{>0}^m$ because $\mathcal{M}(S,\Gamma,\vec{b})$ does not depend upon the geometry of $S$, and we use \eqref{levelset} to describe the pullback Weil-Petersson structure $\omega_{S,\Gamma,\vec{b}}:=\pi_\Gamma^*\omega_{S,\vec{b}}$ on $\mathcal{M}(S,\Gamma,\vec{b})$. Let $\ell_i:\mathcal{M}(S,\Gamma,\vec{b})\rightarrow\mathbb{R}_{>0}$ denote the length of the $i$-th geodesic in $f_*\Gamma$ and let $\tau_i+c_i\mathbb{Z}\in\mathbb{R}/c_i\mathbb{Z}$ denote the twist parameter for $\gamma_i$. Note that $\mathrm{d}\tau_i$ is well-defined. Then, the Weil-Petersson form on $\mathcal{M}(S,\Gamma,\vec{b})$ is given by:
\begin{align}\label{wpform}
\omega_{S,\Gamma,\vec{b}}
=\sum_{j=1}^m\mathrm{d}\ell_j\wedge\mathrm{d}\tau_j+\sum_{i=1}^k\omega_{S,(\vec{c}_i,\vec{b}_i)}.
\end{align}
As with $\mathcal{M}(S,\vec{b})$, we may take the top-exterior product of $\omega_{S,\Gamma,\vec{b}}$ to obtain the Weil-Petersson volume form on $\mathcal{M}(S,\Gamma,\vec{b})$:
\begin{align}\label{wpvolumeform}
\Omega_{S,\Gamma,\vec{b}}=\pi_\Gamma^*\Omega_{S,\vec{b}}
=\bigwedge_{j=1}^m(\mathrm{d}\ell_j\wedge\mathrm{d}\tau_j)\wedge\bigwedge_{i=1}^k \Omega_{S,(\vec{c}_i,\vec{b}_i)}
.
\end{align}
As usual, we omit $\vec{b}$ if it is the (length $n$) zero vector.

\section{McShane identities}\label{mcshanesection}\index{McShane identity}
The following theorem is a rephrasing of McShane's original identity for one-cusped hyperbolic tori \cite{mcshane}:
\begin{theorem}
For any marked one-cusped hyperbolic torus $[X,f]\in\mathcal{T}(S_{1,1})$, let $\mathscr{C}(S_{1,1})$ denote the collection of (non-peripheral) simple closed geodesics on $S_{1,1}$, then
\begin{align}\label{mcshane}
\sum_{\alpha\in\mathscr{C}(S_{1,1})}\frac{2}{1+\exp{\ell_\alpha}([X,f])}=1.
\end{align}
\end{theorem}

\begin{remark}
Each summand in the above series has the following geometric interpretation: the probability that a geodesic launched from the cusp in $X$ will self-intersect before hitting $f_*\alpha$ is precisely $\frac{2}{1+\exp{\ell_\alpha}}$.
\end{remark}

\subsection{Sample application to volume integration}\label{volumesamplesubsec}
Fix a simple closed geodesic $\gamma$ on $S_{1,1}$ and note that $\mathscr{C}(S_{1,1})=\mathrm{Mod}(S_{1,1})\cdot\gamma$. Mirzakhani saw that the $\pi_\gamma$-pushforward of the Weil-Petersson measure for $\Omega_{S_{1,1}}$ weighted by the function $\frac{2}{1+\exp{\ell_\gamma}}:\mathcal{M}(S_{1,1},\gamma)\rightarrow\mathbb{R}_{>0}$ is precisely the Weil-Petersson volume measure induced by $\Omega_{S_{1,1}}$. To see this over a point $[X]\in\mathcal{M}(S_{1,1})$:
\begin{align}\label{unwrap}
(\pi_{\gamma})_*\left(\frac{2\;\Omega_{S_{1,1},\gamma}}{1+\exp{\ell_\gamma}}\right)
&
=\sum_{[X,f_*h_*\gamma]\in\pi_\gamma^{-1}([X])}\frac{2\;\Omega_{S_{1,1}}}{1+\exp{\ell_\gamma([X,f_*h_*\gamma])}}\\
&
=\sum_{\alpha\in\mathrm{Mod}(S_{1,1})\cdot\gamma}\frac{2\;\Omega_{S_{1,1}}}{1+\exp{\ell_\alpha}([X,f_*\alpha])}\notag\\
&
=\sum_{\alpha\in\mathscr{C}(S_{1,1})}\frac{2\;\Omega_{S_{1,1}}}{1+\exp{\ell_\alpha}([X,f])}
=\Omega_{S_{1,1}}.\notag
\end{align}
This in turn means that the following integrals are equivalent:
\begin{align*}
\int_{\mathcal{M}(S_{1,1})}\Omega_{S_{1,1}}=\int_{\mathcal{M}(S_{1,1}\gamma)}\frac{2}{1+\exp{\ell_\gamma}}\;\Omega_{S_{1,1},\gamma}.
\end{align*}
This is our motivation for studying McShane identities: they allow us to unwrap the Weil-Petersson volume $V_{1,1}$ of $\mathcal{M}(S_{1,1})$ as an integral of a function over a topologically simpler moduli space. In this particular case, cutting $S_{1,1}$ along $\gamma$ results in a pair of pants with boundary lengths $(0,\ell_\gamma,\ell_\gamma)$. Hence, \eqref{levelset} tells us that the intermediate moduli space $\mathcal{M}(S_{1,1},\gamma)$ is given by
\begin{align*}
\left\{ (\ell,\tau+\ell\mathbb{Z})\mid \ell\in\mathbb{R}_{>0}\text{ and }\tau+\ell\mathbb{Z}\in\mathbb{R}/\ell\mathbb{Z}\right\},
\end{align*}
with volume form:
\begin{align*}
\Omega_{S_{1,1},\gamma}=\mathrm{d}\ell\wedge\mathrm{d}\tau.
\end{align*}
Therefore, the WP-volume of $\mathcal{M}(S_{1,1})$ \emph{should be}
\begin{align*}
\int_{0}^\infty\frac{2\ell\;\mathrm{d}\ell}{1+\exp\ell}=\frac{\pi^2}{6}.
\end{align*}
However, since one-cusped tori all have an order $2$ isometry called the hyperelliptic involution, we halve $\frac{\pi^2}{6}$ to derive that $V_{1,1}=\frac{\pi^2}{12}$. This is explained more precisely in Remark~\ref{hyperelliptic}.

\subsection{McShane identities for bordered hyperbolic surfaces}
In order to generalize this volume integration strategy, Mirzakhani generalized McShane identities for any bordered hyperbolic surface.

\begin{theorem}[McShane identity]\index{McShane identity}
Let $S$ be a bordered hyperbolic surface with genus $g$ and $n$ labeled geodesic boundaries $\beta_1,\ldots,\beta_n$ of lengths $b_1,\ldots,b_n$, let
\begin{itemize}
\item
$\mathscr{C}_i(S)$ be the collection of simple closed geodesics $\alpha$ which, along with the boundaries $\beta_1$ and $\beta_i$, bound a pair of pants in $S$;
\item
$\mathscr{C}(S)$ be the collection of unordered pairs of simple closed geodesics $\{\alpha_1,\alpha_2\}$ which, along with the boundary $\beta_1$, bound a pairs of pant in $S$.
\end{itemize}
Then, for any marked surface $[X,f]\in\mathcal{T}(S,\vec{b})$, we have the following identity:
\begin{align}\label{mirzakhaniidentity}
b_1=&
\sum_{i=2}^n\sum_{\alpha\in\mathscr{C}_i(S)}
\mathrm{sid}(b_1,b_i;\ell_\alpha([X,f]))
\notag\\
+&
\sum_{\{\alpha_1,\alpha_2\}\in\mathscr{C}(S)} 
\mathrm{mid}(b_1;\ell_{\alpha_1}([X,f]),\ell_{\alpha_2}([X,f])),
\end{align}
where the functions $\mathrm{sid}(\;\cdot\;,\;\cdot\;;\;\cdot\;)$ and $\mathrm{mid}(\;\cdot\;;\;\cdot\;,\;\cdot\;)$ are defined by:
\begin{align}
\mathrm{sid}(b_1,b_i;\ell_\alpha)&:=\log\left(\frac{\cosh\frac{\ell_\alpha}{2}+\cosh\frac{b_1+b_i}{2}}{\cosh\frac{\ell_\alpha}{2}+\cosh\frac{b_1-b_i}{2}}\right)\\
\mathrm{mid}(b_1;\ell_{\alpha_1},\ell_{\alpha_2})&:=2\log\left(\frac{\exp\tfrac{b_1}{2}+\exp\tfrac{\ell_{\alpha_1}+\ell_{\alpha_2}}{2}}{\exp\tfrac{-b_1}{2}+\exp\tfrac{\ell_{\alpha_1}+\ell_{\alpha_2}}{2}}\right).
\end{align}
Note that $\{\alpha_1,\alpha_2\}\in\mathscr{C}(S)$ may contain boundary geodesics.
\end{theorem}

\begin{remark}
We introduce the functions $\mathrm{sid}$ and $\mathrm{mid}$ for minor expositional reasons; they are related to Mirzakhani's $\mathcal{D}$ and $\mathcal{R}$ functions by:
\begin{align}
\mathcal{D}(x,y,z)=\mathrm{mid}(x;y,z)\text{ and }\mathcal{R}(x,y,z)=\mathrm{sid}(x,y;z)+\mathrm{mid}(x,y;z).
\end{align}
\end{remark}

\subsubsection{Proof strategy.}
McShane identities for bordered hyperbolic surfaces may be derived by splitting the length $b_1$ of boundary $\beta_1$ on $[X,f]$ into a countable sum. The idea is to orthogonally shoot out geodesic rays from points on boundary $\beta_1$ and hence partition $\beta_1$ based on the behavior of these orthogeodesic rays. Specifically, starting from a point $x\in\beta_1$, precisely one of three things happens:
\begin{enumerate}
\item
the geodesic ray hits a boundary $\beta_i\neq\beta_1$ without self-intersecting,
\item
the geodesic ray hits $\beta_1$ or self-intersects,
\item
this geodesic ray goes on forever without every self-intersecting.
\end{enumerate}

\noindent\textbf{Case 1:} let $\rho_x$ denote the geodesic arc emanating from $x\in\beta_1$ that hits the $i$-th boundary $\beta_i$. Then, $\beta_1\cup\rho_x\cup\beta_i$ may be fattened up to a unique (homotopy equivalent) pair of pants. The arc $\rho_x$ must be wholly contained in $P_x$, or else forms a hyperbolic $2$-gon --- a geometric impossibility.\\[1em]

\begin{figure}[h]
\begin{center}
\includegraphics[scale=0.25]{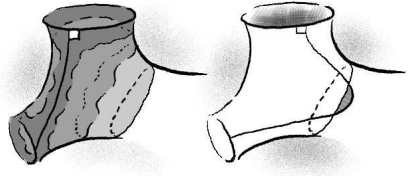}
\end{center}
\caption{From left to right: fattening up $\beta_1\cup\rho_x\cup\beta_i$; an impossible $2$-gon.}
\label{pg11a}
\end{figure}

\noindent\textbf{Case 2:} let $\rho_x$ denote the geodesic arc emanating from $x\in\beta_1$ up to its intersection with $\beta_1$ or its first point of self-intersection. Then, $\beta_1\cup\rho_x$ may be fattened up to a unique (homotopy equivalent) pair of pants. The arc $\rho_x$ must be wholly contained in $P_x$, or else forms either a hyperbolic 2-gon or a hyperbolic triangle with internal angles strictly greater than $\pi$ --- both geometric impossibilities.\\[1em]

\begin{figure}[h]
\begin{center}
\includegraphics[scale=0.25]{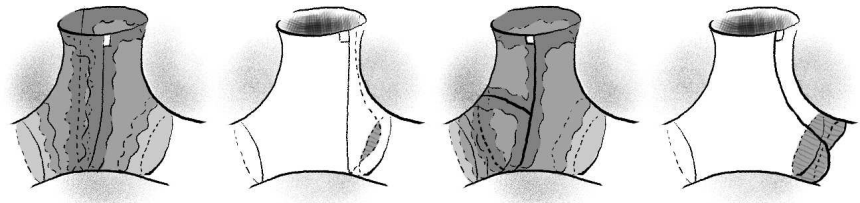}
\end{center}
\caption{From left to right: fattening a border-hitting $\beta_1\cup\rho_x$; an impossible $2$-gon; fattening self-intersecting $\beta_1\cup\rho_x$; an impossible triangle (cut the annulus).}
\label{pg11bc}
\end{figure}

\noindent\textbf{Case 3:} there are uncountably many such simple orthogeodesic rays. However, the Birman-Series geodesic sparsity theorem\index{Birman-Series theorem} \cite{birmanseries} tells us that they occupy a set of measure $0$ on $\beta_1$. More accurately, the Birman-Series theorem says that the collection of simple geodesics on a closed hyperbolic surface has Hausdorff dimension $1$. Double a bordered hyperbolic surface $X$ by gluing an isometric (but orientation-reversed) copy $\bar{X}$ to $X$ along correspondingly labeled borders. Orthogeodesics rays on $X$ then glue to corresponding orthogeodesics rays on $\bar{X}$ to give simple bi-infinite geodesics on the double. The Birman-Series theorem then asserts that the set of points occupied by these bi-infinite geodesics has measure $0$ on $X$. This in turn means that the restriction of these simple bi-infinite geodesics to a collar neighborhood of boundary $\beta_1$ (as a subset of the double) occupies $0$ area. However, this collar neighborhood has the structure of an interval times $\beta_1$, and so we see that the set of points on boundary $1$ that launch simple orthogeodesic rays has measure $0$ with respect to the length measure on the boundary.

We conclude therefore that almost every $x\in\beta_1$ belongs either to Case 1 or Case 2, and in these two cases, the geodesic arc $\rho_x$ lies on a unique pair of pants $P$ in $X$. This gives us a natural decomposition of the total measure $b_1$ of $\beta_1$ as an infinite sum over pairs of pants $P$ embedded in $X$.

\subsubsection{Orthogeodesics on pairs of pants.}\index{pair of pants!orthogeodesics}

\begin{figure}
\begin{center}
\includegraphics[scale=0.125]{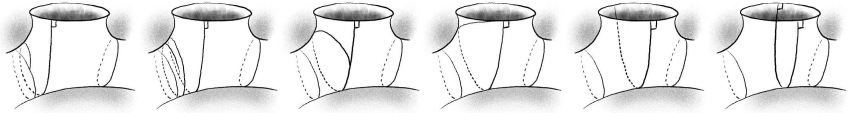}
\end{center}
\caption{A ``movie'' of various types of orthogeodesic behavior.}
\label{pg12}
\end{figure}

Given a pair of pants $P\subset X$ with boundaries $\beta_1,\alpha_1,\alpha_2$, there are precisely four simple infinite orthogeodesic rays contained in $P$ (Figure~\ref{pg12}, second from left). The end points of these four rays partition $\beta_1$ into four intervals.
Geodesic rays launched from the interval closest to $\alpha_1$ necessarily hit $\alpha_1$ (before possibly self-intersecting; Figure~\ref{pg12}, leftmost). The width of this \emph{side interval} (computed with a little hyperbolic trigonometry \cite[Thm.~2.3.1]{buser}) is
\begin{align}\label{sidegap}
\mathrm{sid}(b_1,\ell_{\alpha_1};\ell_{\alpha_2})=
\log\left(\frac{\cosh\frac{\ell_{\alpha_2}}{2}+\cosh\frac{b_1+\ell_{\alpha_1}}{2}}{\cosh\frac{\ell_{\alpha_2}}{2}+\cosh\frac{b_1-\ell_{\alpha_1}}{2}}\right).
\end{align}
Terms of the form \eqref{sidegap} arise in Case 1, that is: for embedded pairs of pants $P\subset X$ which contain two distinct boundary geodesics $\beta_1,\beta_i$. Since pairs of pants $P\subset X$ may be given by specifying their boundary geodesics in $X$, the collection of pairs of pants $P$ bordered by $\beta_1$ and $\beta_i$ precisely corresponds to $\mathscr{C}_i(S)$. Replacing $\ell_{\alpha_1}$ with $b_i$ and $\ell_{\alpha_2}$ with $\ell_\alpha$ in \eqref{sidegap} produces the correct summand $\mathrm{sid}(b_1,b_i;\ell_\alpha)$ in Mirzakhani's identity.

By symmetry, the above statements also hold for the interval closest to $\alpha_2$ upon switching the roles of $\alpha_1$ and $\alpha_2$.

Orthogeodesic rays launched from the two remaining intervals either self-intersect (before possibly leaving $P$; Figure~\ref{pg12}, center-left) or hit $\beta_1$ (Figure~\ref{pg12}, right). The width of each of these two \emph{middle intervals} is:
\begin{align}\label{innergap}
\tfrac{1}{2}\mathrm{mid}(b_1;\ell_{\alpha_1},\ell_{\alpha_2})=\log\left(\frac{\exp\tfrac{b_1}{2}+\exp\tfrac{\ell_{\alpha_1}+\ell_{\alpha_2}}{2}}{\exp\tfrac{-b_1}{2}+\exp\tfrac{\ell_{\alpha_1}+\ell_{\alpha_2}}{2}}\right).
\end{align}
Terms of the form \eqref{innergap} arise in Case 2, that is: for embedded pairs of pants $P\subset X$ which contain $\beta_1$. This corresponds to $\mathscr{C}(S)$. Doubling \eqref{innergap} due to there being two such intervals for $P$ produces the correct summand $\mathrm{mid}(b_1;\ell_{\alpha_1},\ell_{\alpha_2})$ in the bordered McShane identity. 

This completes the proof of McShane identity for bordered hyperbolic surfaces.

\subsubsection{Simple infinite orthogeodesic rays.}

As a minor aside, we point out that Mirzakhani also gives a detailed analysis of what occurs in Case 3.

\begin{theorem}[{\cite[Thm.~4.5, 4.6]{mirz1}}]\label{case3thm}
The set of points on $x\in\beta_1$ corresponding to simple infinite orthogeodesics is homeomorphic to the Cantor set union countably many isolated points. Specifically, if the orthogeodesic ray $\rho_x$ emanating from $x$
\begin{enumerate}
\item spirals to a boundary curve $\beta_{i\neq1}$, then $x$ is an isolated point;
\item spirals to a simple closed geodesic in the interior of $X$, then $x$ is a boundary point of the Cantor set;
\item does not spiral to a simple closed curve, then $x$ is a non-boundary point of the Cantor set.
\end{enumerate}
\end{theorem}
Mirzakhani's proof of the above result is slightly technical, but the result itself is geometrically unsurprising. In the course of establishing her McShane identity, we have seen that the Case 1 and Case 2 points $x\in\beta_1$ each lie within precisely one of four (open) intervals on a certain pair of pants $\{\alpha_1,\alpha_2\}\in\mathscr{C}(S)$. Case 3 points are the points on $\beta_1$ that still remain after removing these open intervals.

When a pair of pants $\{\alpha_1,\alpha_2\}$ contains a boundary geodesics $\alpha_1=\beta_i$ distinct from $\beta_1$, we need to remove three of the intervals --- leaving the interval closest to $\alpha_2$. The end points of the two simple infinite orthogeodesic rays wedged in between these three intervals are obviously isolated points, and these rays spiral to $\beta_i$. All isolated points arise in this way.

If we add these isolated points to the points that we remove from $\beta_1$, then for every pair of pants $P\in \mathscr{C}_i(S)$, we remove one long open interval (containing three of the original intervals) and for every pair of pants $P\in\mathscr{C}(S)-\cup\mathscr{C}_i(S)$, we remove two open intervals. This process of removing intervals from $\beta_1$ is essentially akin to how the standard Cantor set is constructed, and it should be expected that the remnant collection of points is a Cantor set. Moreover, this description tells us that the boundary points of the intervals that we remove correspond to orthogeodesic rays which spiral to an interior simple closed geodesic --- as was asserted in statement (2).

Note that the existence of the Case 3 Cantor set is one reason for which we needed to invoke the Birman-Series theorem. After all, the measure of a Cantor set on $\beta_1$ can take any value in $[0,b_1]$.

\section{Weil-Petersson volume computation}\label{volumesec}

Let $V_{g,n}(\vec{b})$ denote the Weil-Petersson volume of the moduli space $\mathcal{M}(S_{g,n},\vec{b})$.

\subsection{Mirzakhani's volume recursion formula}

\subsubsection{Derivation and $V_{1,3}(b_1,b_2,b_3)$.}
We summarize the key steps of Mirzakhani's volume computation procedure, while giving a step-by-step calculation of $V_{1,3}(\vec{b})$ as an illustrative example.\\[1.5em]
\noindent\textbf{Step 1:} Rearrange the McShane identity into mapping class group orbits of \emph{ordered} tuples of curves. This prepares the McShane identity in a form conducive to the integral unwrapping we saw in \eqref{unwrap}.

First identify elements of the form $\{\beta_i,\alpha\}\in\mathscr{C}(S)$ and $\alpha\in\mathscr{C}_i(S)=\mathrm{Mod}(S)\cdot\gamma_i$ since they both correspond to the pair of pants on $S$ bordered by $\{\beta_1,\beta_i,\alpha\}$. Gather the summands of the form $\mathrm{mid}(b_1;b_i,\ell_\alpha)$ over $\{\beta_i,\alpha\}$ with the corresponding summand $\mathrm{sid}(b_1,b_i;\ell_\alpha)$ over $\alpha\in\mathscr{C}_i(S)$ to get series of the form
\begin{align*}
\sum_{i=2}^n\sum_{\alpha\in\mathrm{Mod}(S)\cdot\gamma_i}[\mathrm{sid}(b_1,b_i;\ell_\alpha)+\mathrm{mid}(b_1;b_i,\ell_\alpha)].
\end{align*}

The remaining elements of $\mathscr{C}(S)-\cup\mathscr{C}_i(S)$ correspond to pairs of pants whose non-$\beta_1$ borders are on the interior of $S$. We break up the summands $\mathrm{mid}(b_1;\ell_{\alpha_1},\ell_{\alpha_2})$ over $\{\alpha_1,\alpha_2\}\in\mathscr{C}(S)-\cup\mathscr{C}_i(S)$ into two summands of the form $\frac{1}{2}\mathrm{mid}(b_1;\ell_{\alpha_1},\ell_{\alpha_2})$ over $(\alpha_1,\alpha_2)$ and $(\alpha_2,\alpha_1)$, thereby enabling us to sum over
\begin{align*}
\{(\alpha_1,\alpha_2)\mid \{\alpha_1,\alpha_2\}\in\mathscr{C}(S)-\cup\mathscr{C}_i(S)\}.
\end{align*}
Partition the new summation index set of ordered geodesic pairs into mapping class group orbits, and gather the summands accordingly. 

\begin{remark}
To determine whether two ordered curves $(\alpha_1,\alpha_2)$ and $(\alpha_1',\alpha_2')$ are in the same mapping class group orbit, check if $(\alpha_1,\alpha_2)$ and $(\alpha_1',\alpha_2')$ respectively decompose $S$ into surfaces with topologically equivalent connected components with matching boundary labels and where $\alpha_i$ matches with $\alpha_i'$.
\end{remark}

\begin{figure}
\begin{center}
\includegraphics[scale=0.3]{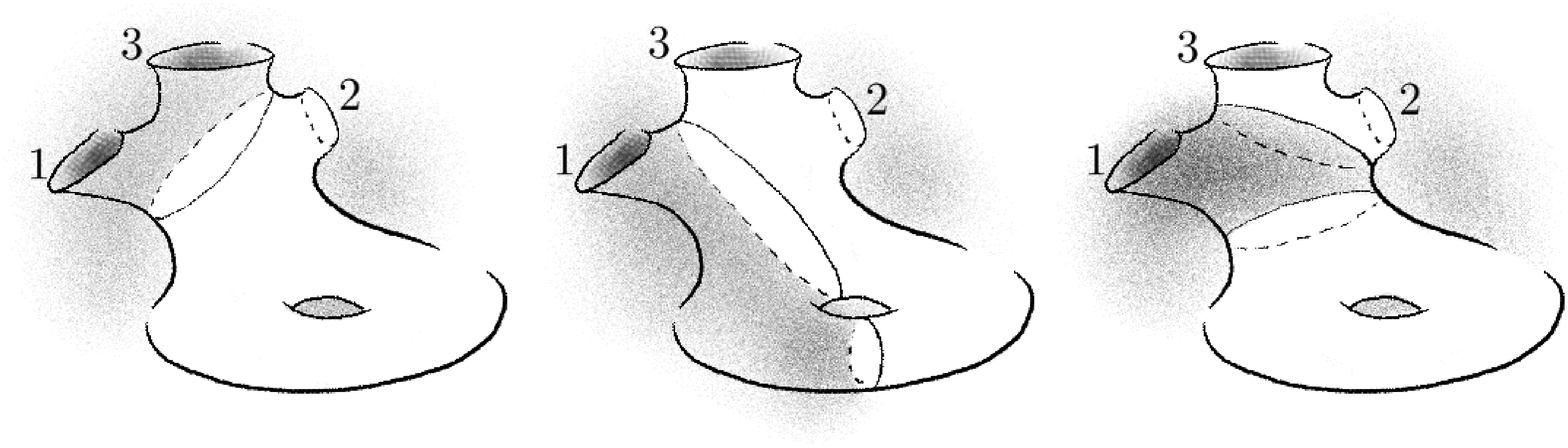}
\end{center}
\caption{Three topologically distinct ways of excising a pair of pants containing boundary $\beta_1$ from a thrice-holed torus.}
\label{pg15}
\end{figure}

\begin{example}
Consider a thrice-holed hyperbolic torus $S=S_{1,3}$ with boundary lengths $\vec{b}=(b_1,b_2,b_3)\in\mathbb{R}_{>0}^3$, and recall that Mirzakhani's McShane identity \eqref{mirzakhaniidentity} for marked surfaces $[X,f]\in\mathcal{T}(S,\vec{b})$ is a series summed over certain sets $\mathscr{C}_2(S),\mathscr{C}_3(S)$ and $\mathscr{C}(S)$. Choose arbitrary elements $\gamma_i\in\mathscr{C}_i(S)$, then $\mathscr{C}_i(S)=\mathrm{Mod}(S)\cdot\gamma_i$.

The ordered summation index
\begin{align*}
\{(\alpha_1,\alpha_2)\mid \{\alpha_1,\alpha_2\}\in\mathscr{C}(S)-(\mathcal(S)_2(S)\cup\mathscr{C}_3(S))\}
\end{align*}
partitions into three mapping class orbits:
\begin{align*}
\mathrm{Mod}\cdot(\gamma^{\mathrm{con}},\delta^{\mathrm{con}})\sqcup\mathrm{Mod}(S)\cdot(\gamma^{\mathrm{dcon}},\delta^{\mathrm{dcon}}) \sqcup \mathrm{Mod}(S)\cdot(\delta^{\mathrm{dcon}},\gamma^{\mathrm{dcon}}),
\end{align*}
where $\gamma^{\mathrm{con}},\delta^{\mathrm{con}}$ are chosen so that excising the pairs of pants bordered by $\{\beta_1,\gamma^{\mathrm{con}},\delta^{\mathrm{con}}\}$ from $S$ results in a connected surface (Figure~\ref{pg15} center), and where $\gamma^{\mathrm{con}},\delta^{\mathrm{con}}$ are chosen so that excising the pairs of pants bordered by $\{\beta_1,\gamma^{\mathrm{con}},\delta^{\mathrm{con}}\}$ results in a disconnected surface (Figure~\ref{pg15} right).
Thus, the McShane identity may be rearranged as follows:
\begin{align}\label{13mcshane}
b_1=&\sum_{i=2}^3\sum_{\alpha\in\mathrm{Mod}(S)\cdot\gamma_i}[\mathrm{sid}(b_1,b_i;\ell_\alpha)+\mathrm{mid}(b_1;b_i,\ell_\alpha)]\\
&+\frac{1}{2}\sum_{\{\alpha_1,\alpha_2\}\in\mathrm{Mod}(S)\cdot\{\gamma^{\mathrm{con}},\delta^{\mathrm{con}}\}}\mathrm{mid}(b_1;\ell_{\alpha_1},\ell_{\alpha_2})\notag\\
&+\frac{1}{2}\sum_{(\alpha_1,\alpha_2)\in\mathrm{Mod}(S)\cdot(\gamma^{\mathrm{dcon}},\delta^{\mathrm{dcon}})}\mathrm{mid}(b_1;\ell_{\alpha_1},\ell_{\alpha_2})\notag\\
&+\frac{1}{2}\sum_{(\alpha_1,\alpha_2)\in\mathrm{Mod}(S)\cdot(\delta^{\mathrm{dcon}},\gamma^{\mathrm{dcon}})}\mathrm{mid}(b_1;\ell_{\alpha_1},\ell_{\alpha_2}).\notag
\end{align}
\end{example}

\noindent\textbf{Step 2:} Integrate the McShane identity over moduli space. The left hand side is $b_1 V_{g,n}(\vec{b})$ and the right hand side \emph{unwraps} as integrals over various intermediate moduli spaces. Specifically, we use the fact that the pushforward of the weighted WP volume measure $[\mathrm{sid}(b_1,b_2;\ell_{\gamma_i})+\mathrm{mid}(b_1;b_i,\ell_{\gamma_i})]\cdot\Omega_{S,\gamma_i,\vec{b}}$ onto $\mathcal{M}(S,\vec{b})$ with respect to the forgetful map $\pi_{\gamma_i}$ is precisely given by
\begin{align*}
\left(\sum_{\alpha\in\mathrm{Mod}(S)\cdot\gamma_i}[\mathrm{sid}(b_1,b_2;\ell_{\alpha_i})+\mathrm{mid}(b_1;b_i,\ell_{\alpha_i})]\right)\cdot\Omega_{S,\vec{b}}
\end{align*}
to unwrap to $\mathcal{M}(S,\gamma_i,\vec{b})$. And we use the fact that the pushforward measure of the weighted WP volume measure $\frac{1}{2}\mathrm{mid}(b_1;\ell_{\gamma},\ell_{\delta})\cdot\Omega_{S,(\gamma,\delta),\vec{b}}$ onto $\mathcal{M}(S,\vec{b})$ with respect to the forgetful map $\pi_{(\gamma,\delta)}$ is
\begin{align*}
\frac{1}{2}\left(\sum_{(\alpha_1,\alpha_2)\in\mathrm{Mod}(S)\cdot(\gamma,\delta)}\mathrm{mid}(b_1;\ell_{\alpha_1},\ell_{\alpha_2})\right)\cdot\Omega_{S,\vec{b}}
\end{align*}
to unwrap to $\mathcal{M}(S,(\gamma,\delta),\vec{b})$.

\begin{example}
Integrating \eqref{13mcshane} over $\mathcal{M}(S,\vec{b})$, we obtain
\begin{align}\label{13rough}
b_1\cdot V_{1,3}(\vec{b})=&\sum_{i=2}^3\int_{\mathcal{M}(S,\gamma_i,\vec{b})}[\mathrm{sid}(b_1,b_i;\ell_\alpha)+\mathrm{mid}(b_1;b_i,\ell_\alpha)]\;\Omega_{S,\gamma_i,\vec{b}}\\
&+\frac{1}{2}\int_{\mathcal{M}(S,(\gamma^{\mathrm{con}},\delta^{\mathrm{con}}),\vec{b})}\mathrm{mid}(b_1;\ell_{\gamma^{\mathrm{con}}},\ell_{\delta^{\mathrm{con}}})\;\Omega_{S,(\gamma^{\mathrm{con}},\delta^{\mathrm{con}}),\vec{b}}\notag\\
&+\frac{1}{2}\int_{\mathcal{M}(S,(\gamma^{\mathrm{dcon}},\delta^{\mathrm{dcon}}),\vec{b})}\mathrm{mid}(b_1;\ell_{\gamma^{\mathrm{dcon}}},\ell_{\delta^{\mathrm{dcon}}})\;\Omega_{S,(\gamma^{\mathrm{dcon}},\delta^{\mathrm{dcon}}),\vec{b}}\notag\\
&+\frac{1}{2}\int_{\mathcal{M}(S,(\delta^{\mathrm{dcon}},\gamma^{\mathrm{dcon}}),\vec{b})}\mathrm{mid}(b_1;\ell_{\delta^{\mathrm{dcon}}},\ell_{\gamma^{\mathrm{dcon}}})\;\Omega_{S,(\delta^{\mathrm{dcon}},\gamma^{\mathrm{dcon}}),\vec{b}}\notag.
\end{align}
\end{example}

\noindent\textbf{Step 3:} We gather the intermediate moduli spaces we encounter in Step 2 into three groups. For each of these groups, we analyze the structure of an arbitrary level set $\ell_\Gamma (\vec{x})\subset\mathcal{M}(S,\Gamma,\vec{b})$ with equation~\eqref{levelset} and use \eqref{wpvolumeform} to express the volume of $\mathcal{M}(S,\vec{b})$ in terms of integrals of functions built from the volumes of lower dimensional moduli spaces of bordered surfaces. This is the \emph{recursive} aspect of Mirzakhani's formula.

The first group is made up of intermediate moduli spaces of the form $\mathcal{M}(S,\gamma_i,\vec{b})$ --- these moduli spaces correspond to the summands over $\mathscr{C}_i(S)=\mathrm{Mod}(S)\cdot\gamma_i$. Excising the pair of pants bordered by $\{\beta_1,\beta_i,\gamma\}$ from $S$ leaves a (connected) genus $g$ with $n-1$ borders $\{\gamma,\beta_2,\ldots,\beta_n\}$ There are $n-1$ integrals of this type, one for each boundary index $i=2,\ldots, n$.

\begin{example}
Invoking \eqref{levelset}, we see that the level set $\ell_{\gamma_2}^{-1}(x)\subset\mathcal{M}(S,\gamma_2,\vec{b})$ is
\begin{align*}
\mathbb{R}/x\mathbb{Z}\times\mathcal{M}(S_{0,3},(x,b_1,b_2))\times\mathcal{M}(S_{1,2},(x,b_3)).
\end{align*}
But $\mathcal{M}(S_{0,3},(x,b_1,b_2))$ is just a single point, and is assigned volume $1$. The $i=2$ integral in the top line of \eqref{13rough} therefore transforms to:
\begin{align}
\int_0^\infty x\cdot V_{1,2}(x,b_3) \cdot [\mathrm{sid}(b_1,b_2;x)+\mathrm{mid}(b_1;b_2,x)]\;\mathrm{d}x.
\end{align}
The $i=3$ integral on the top line transforms alike by symmetry.
\end{example}

The second group is made up of intermediate moduli spaces of the form $\mathcal{M}(S,(\gamma^{\mathrm{con}},\delta^{\mathrm{con}}),\vec{b})$, where excising the hyperbolic pair of pants on $S$ bordered by $\{\beta_1,\gamma^{\mathrm{con}},\delta^{\mathrm{con}}\}$ results in a genus $g-1$ connected surface bordered by $n+1$ boundaries $\{\gamma^{\mathrm{con}},\delta^{\mathrm{con}},\beta_2,\ldots,\beta_n\}$. There is precisely $1$ integral of this type.

\begin{example}
The level set $\ell_{(\gamma^{\mathrm{con}},\delta^{\mathrm{con}})}^{-1}(x,y)$ in the intermediate moduli space $\mathcal{M}(S,(\gamma^{\mathrm{con}},\delta^{\mathrm{con}}),\vec{b})$ may be identified with
\begin{align*}
\mathbb{R}/x\mathbb{Z}\times\mathbb{R}/y\mathbb{Z}\times\mathcal{M}(S_{0,3},(x,y,b_1))\times\mathcal{M}(S_{0,4}(x,y,b_2,b_3)).
\end{align*}
Thus, the second line integral transforms to
\begin{align}
\frac{1}{2}\int_0^\infty\int_0^\infty xy\cdot V_{0,4}(x,y,b_2,b_3)\cdot \mathrm{mid}(b_1;x,y)\;\mathrm{d}x\;\mathrm{d}y.
\end{align}
\end{example}

The third and final group consists of intermediate moduli spaces of the form $\mathcal{M}(S,\gamma^{\mathrm{dcon}},\delta^{\mathrm{dcon}}),\vec{b})$, where excising the pair of pants on $S$ bordered by $\{\beta_1,\gamma^{\mathrm{dcon}},\delta^{\mathrm{dcon}}\}$ results in two connected components $S_1,S_2$. Denote the respective genus of these two surfaces $S_1,S_2$ by $g_1,g_2$ and let their respective boundary components by $\{\gamma^{\mathrm{dcon}}\}\cup\{\beta_i\}_{i\in I_1}$ and $\{\delta^{\mathrm{dcon}}\}\cup\{\beta_j\}_{j\in I_2}$, we know immediately that $g_1+g_2=g$, the union $I_1\sqcup I_2$ of these disjoint index sets is equal to $\{2,\ldots,n\}$ and $2g_1-2+|I_1|+1, 2g_2-2+|I_2|+1\geq 0$. Conversely, given $g_1,g_2\geq0$ and disjoint (potentially empty) index sets $I_1,I_2$ satisfying the above conditions, there are geodesics $\{\gamma^{\mathrm{dcon}},\delta^{\mathrm{dcon}}\}$ on $S$ so that excising the pair of pants bordered by $\{\beta_1,\gamma^{\mathrm{dcon}},\delta^{\mathrm{dcon}}\}$ leaves two boundary-labeled surfaces respectively homeomorphic to $S_1$ and $S_2$.

\begin{example}
The level set $\ell_{(\gamma^{\mathrm{dcon}},\delta^{\mathrm{dcon}})}^{-1}(x,y)\subset\mathcal{M}(S,(\gamma^{\mathrm{dcon}},\delta^{\mathrm{dcon}}),\vec{b})$ is
\begin{align*}
\mathbb{R}/x\mathbb{Z}\times\mathbb{R}/y\mathbb{Z}\times\mathcal{M}(S_{0,3},(x,y,b_1))\times\mathcal{M}(S_{0,3}(x,b_2,b_3))\times\mathcal{M}(S_{1,1}(y)).
\end{align*}
Therefore, the last two lines in \eqref{13rough} respectively become
\begin{align}
\frac{1}{2}\int_0^\infty\int_0^\infty xy\cdot V_{0,3}(x,b_2,b_3)\cdot V_{1,1}(y)\cdot \mathrm{mid}(b_1;x,y)\;\mathrm{d}x\;\mathrm{d}y
\end{align}
and the same integral (with $x$ and $y$ switched).
\end{example}

\noindent\textbf{Step 4:} To complete the volume calculation, we differentiate $b_1V_{g,n}(\vec{b})$ with respect to $b_1$. This reduces $\mathrm{sid}(b_1,b_i;x)$ and $\mathrm{mid}(b_1;x,y)$ into functions resembling a function we previously integrated in Subsection~\ref{volumesamplesubsec} to obtain the volume of $\mathcal{M}(S_{1,1})$. Define the function
\begin{align*}
H(s,t):=\frac{1}{1+\exp\frac{s+t}{2}}+\frac{1}{1+\exp\frac{s-t}{2}}.
\end{align*}
Differentiating $\mathrm{mid}(b_1;x,y)$ with respect to $b_1$ gives
$H(x+y,b_1)$
and differentiating $\mathrm{mid}(b_1;b_i,x)+\mathrm{sid}(b_1,b_i;x)$ with respect to $b_1$ gives
$\tfrac{1}{2}[H(x,b_1+b_i)+H(x,b_1-b_i)]
$.
\begin{example}
After differentiating by $b_1$, the volume $b_1 V_{1,3}(\vec{b})$ becomes:
\begin{align}\label{13fine}
\frac{\partial [2b_1 V_{1,3}(\vec{b})]}{\partial\; b_1}=
&\int_0^\infty x V_{1,2}(x,b_3)\left[H(x,b_1+b_2)+H(x,b_1-b_2)\right]\mathrm{d}x \\
+&\int_0^\infty x V_{1,2}(x,b_2)\left[H(x,b_1+b_3)+H(x,b_1-b_3)\right]\mathrm{d}x \notag\\
+&\int_0^\infty\int_0^\infty xy V_{0,4}(x,y,b_2,b_3) H(x+y,b_1)\; \mathrm{d}x\;\mathrm{d}y \notag\\
+&\int_0^\infty\int_0^\infty xy V_{0,3}(x,b_2,b_3) V_{1,1}(y) H(x+y,b_1)\;\mathrm{d}x\;\mathrm{d}y\notag\\
+&\int_0^\infty\int_0^\infty xy V_{1,1}(x) V_{0,3}(y,b_2,b_3) H(x+y,b_1)\;\mathrm{d}x\;\mathrm{d}y,\notag
\end{align}
where the recursion volume polynomials are
\begin{align*}
V_{1,1}(x_1)&=\tfrac{1}{12}\pi^2+\tfrac{1}{48}x_1^2\\
V_{0,4}(x_1,x_2,x_3,x_4)&=2\pi^2+\tfrac{1}{2}(x_1^2+x_2^2+x_3^2+x_4^2)\\
V_{1,2}(x_1,x_2)&=\tfrac{1}{4}\pi^4+\tfrac{1}{12}(\pi^2b_1^2+\pi^2b_2^2)+\tfrac{1}{96}x_1^2x_2^2+\tfrac{1}{192}(x_1^4+x_2^2),
\end{align*}
and where a factor of $\frac{1}{2}$ has been multiplied out onto the left hand side.
\end{example}

We discuss in Subsection~\ref{polynomialsubsec} and demonstrate in Appendix~\ref{appendixsec} how to explicitly compute these integrals, and show that $\frac{\partial [b_1V_{g,n}(\vec{b})]}{\partial\; b_1}$ is always a polynomial. This allow us to obtain $V_{g,n}(\vec{b})$ up to the addition of some function in $b_2,\ldots,b_n$. This ambiguity is easily resolved by using the fact that $V_{g,n}(\vec{b})$ is symmetric in the $b_i$.

\begin{example}
Integrating out \eqref{13fine}, we get
\begin{align*}
\tfrac{\partial [2b_1 V_{1,3}(\vec{b})]}{\partial\; b_1}=&\tfrac{28}{9}\pi^6+\tfrac{13}{12}(3\pi^4b_1^2+\pi^4b_2^2+\pi^4b_3^2)\\
&+\tfrac{1}{4}(3\pi^2b_1^2b_2^2+3\pi^2b_1^2b_3+\pi^2b_2^2b_3^2)\\
&+\tfrac{1}{12}(5\pi^2b_1^4+\pi^2b_2^4+\pi^2b_3^4)+\tfrac{1}{16}b_1^2b_2^2b_3^2\\
&+\tfrac{1}{96}(5b_1^4b_2^2+3b_1^2b_2^4+5b_1^4b_3^2+3b_1^2b_3^4+b_2^4b_3^2+b_2^2b_3^4)\\
&+\tfrac{1}{576}(7b_1^6+b_2^6+b_3^6).
\end{align*}
Therefore, we obtain the volume for the moduli space of thrice-holed tori:
\begin{align}
V_{1,3}(\vec{b})=&\tfrac{14}{9}\pi^6+\tfrac{13}{24}(\pi^4b_1^2+\pi^4b_2^2+\pi^4b_3^2)\\
&+\tfrac{1}{8}(\pi^2b_1^2b_2^2+\pi^2b_1^2b_3+\pi^2b_2^2b_3^2)\notag\\
&+\tfrac{1}{24}(\pi^2b_1^4+\pi^2b_2^4+\pi^2b_3^4)+\tfrac{1}{96}b_1^2b_2^2b_3^2\notag\\
&+\tfrac{1}{192}(b_1^4b_2^2+b_1^2b_2^4+b_1^4b_3^2+b_1^2b_3^4+b_2^4b_3^2+b_2^2b_3^4)\notag\\
&+\tfrac{1}{1152}(b_1^6+b_2^6+b_3^6).\notag
\end{align}
\end{example}

\subsubsection{The formula.}
Given a set of indices $I\subset\{1,\ldots,n\}$, we use $\hat{\vec{b}}_{I}$ to denote the vector obtained by removing, for every $i\in I$, the $b_i$ term in $\vec{b}$. Then, Mirzakhani's volume recursion formula is:
\begin{align}\label{recursionformula}\index{Weil-Petersson!volume recursion formula}
\frac{\partial [2b_1V_{g,n}(\vec{b})]}{\partial b_1}=\sum_{i=2}^n&\int_0^\infty x V_{g,n-1}(x,\hat{\vec{b}}_{1,i})\left[H(x,b_1+b_i)+H(x,b_1-b_i)\right]\mathrm{d}x \notag\\
+\int_0^\infty&\int_0^\infty xy V_{g-1,n+2}(x,y,\hat{\vec{b}}_1) H(x+y,b_1)\; \mathrm{d}x\;\mathrm{d}y\\
+\sum_{\substack{((g_1,I_1),(g_2,I_2))\\g_1+g_2=g\\I_1\sqcup I_2=\{2,\ldots,n\}}}\int_0^\infty&\int_0^\infty xy V_{g_1,|I_1|+1}(x,\hat{\vec{b}}_{I_2}) V_{g_2,|I_2|+1}(y,\hat{\vec{b}}_{I_1}) H(x+y,b_1)\;\mathrm{d}x\;\mathrm{d}y,\notag
\end{align}
where $g_1,g_2$ are non-negative integers, the $I_1,I_2$ are (possible empty) disjoint sets and $0\leq 2g_1-2+|I_1|,2g_2-2+|I_2|$.
The base cases for the recursion are:
\begin{align*}
V_{0,3}(b_1,b_2,b_3)=1\text{ and }V_{1,1}(b)=\tfrac{\pi^2}{12}+\tfrac{b^2}{48}.
\end{align*}

\begin{figure}[h]
\begin{center}
\includegraphics[scale=.2]{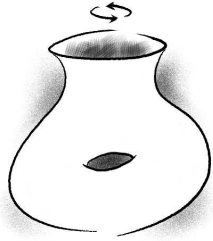}
\end{center}
\caption{Every one-cusped/one-bordered tori has a hyperelliptic involution.}
\label{pg10}
\end{figure}

\begin{remark}\label{hyperelliptic}
There are two conventions for the volume $V_{1,1}(b_1)$ of the moduli space of one-cusped/holed tori. While N\"a\"at\"anen-Nakanishi \cite{nn1,nn2}, Mirzakhani \cite{mirz1} and Wolpert \cite{wolpert1, wolpert2} choose the convention of $\frac{\pi^2}{6}$, we side with Zograf \cite{zograf2} and choose to halve the volume of $\mathcal{M}(S_{1,1},b_1)$ as it simplifies Mirzakhani's recursion formula. Note that this issue applies only to $V_{1,1}(b_1)$. Indeed, $\mathcal{M}(S_{1,1},b_1)$ is the only moduli space of \emph{boundary-labeled} surfaces where every surface has an order $2$ orientation preserving isometry called the \emph{hyperelliptic involution}\index{hyperelliptic involution} (Figure~\ref{pg10}). For every other moduli space, the locus in $\mathcal{M}(S_{g,n},\vec{b})|_{(g,n)\neq(1,1)}$ consisting of boundary-labeled hyperbolic surfaces with hyperelliptic involutions has measure $0$. 
\end{remark}

\subsection{Volume polynomial.}\label{polynomialsubsec}

We prove that Weil-Petersson volumes $V_{g,n}(\vec{b})$ are polynomials, and discuss how integrals of the form \eqref{recursionformula} may be explicitly computed.
\begin{theorem}\label{polynomial}
The Weil-Petersson volume $V_{g,n}(\vec{b})$ of $\mathcal{M}(S_{g,n},\vec{b})$ is a polynomial in $b_1^2,\ldots,b_n^2$ of the following form:
\begin{align}
V_{g,n}(\vec{b})=\sum_{\vec{a}:|\vec{a}|\leq 3g-3+n} C_\vec{a}\cdot b^{2\vec{a}},
\end{align}
where the exponents $\vec{a}=(a_1,\ldots,a_n)$ vary over $\mathbb{Z}_{\geq0}^m$, $b^{2\vec{a}}:=b_1^{2a_1}\times\cdots\times b_n^{2a_n}$ and $C_\vec{a}\in\pi^{6g-6+2n-2|\vec{a}|}\cdot\mathbb{Q}_{>0}$.
\end{theorem}
Mirzakhani proves Theorem~\ref{polynomial} by induction on the number of pairs of pants needed to glue to form a surface $S_{g,n}$, which is also the absolute value of the Euler characteristic $|\chi(S_{g,n})|=2g-2+n$. The two base cases $(g,n)=(0,3),(1,1)$ are true. Assume therefore that moduli spaces for all bordered surfaces $S$ with $|\chi(S)|<2g-2+n$ satisfy Theorem~\ref{polynomial}. We would like to show that 
\begin{align*}
F_{2k+1}(b):=\int_0^\infty x^{2k+1}\cdot H(x,b)\;\mathrm{d}x\text{ and }\int_0^\infty\int_0^\infty x^{2i+1}y^{2j+1}\cdot H(x+y,b)\;\mathrm{d}x\;\mathrm{d}y
\end{align*}
are polynomials in $\mathbb{Q}_{>0}[\pi^2,t^2]$ of the correct degree and with correctly matching powers of $\pi^2$ for each monomial. Furthermore, since 
\begin{align}
\int_0^\infty\int_0^\infty x^{2i+1}y^{2j+1}\cdot H(x+y,b)\;\mathrm{d}x\;\mathrm{d}y=\frac{(2i+1)!(2j+1)!}{(2i+2j+3)!}F_{2i+2j+3}(b),
\end{align}
we only need to consider $F_{2k+1}(b)$. We show in Appendix~\ref{appendixsec} that $F_{2k+1}(b)$ is equal to
\begin{align}
F_{2k+1}(b)={(2k+1)!}\sum_{i=0}^{k+1}\zeta(2i)(2^{2i+1}-4)\frac{b^{2k+2-2i}}{(2k+2-2i)!},
\end{align}
where $\zeta$ denotes the Riemann zeta function. Since $\zeta(0)=\frac{-1}{2}$, the highest order term in this polynomial is positive. More generally, the recursion relation
\begin{align}
\zeta(2i)=\frac{2}{2i+1}\sum_{j=1}^{i-1}\zeta(2j)\zeta(2i-2j)
\end{align} 
and the fact that $\zeta(2)=\frac{\pi^2}{6}$ ensures that the coefficients in $F_{2k+1}(b)$ are of the right form.

\section{Applications}\label{applicationsec}
We discuss here two of Mirzakhani's applications of her volume recursion formula: a proof of Witten's conjecture (Kontsevich's theorem) and her results on the polynomial growth rate of simple closed geodesics on hyperbolic surfaces. Our aim is to give a simple and rough outline of the key ideas and steps in these two applications.

\subsection{Witten's conjecture}
In \cite{witten}, Witten considered two approaches to integrating over the infinite-dimensional space of metrics on a genus $g$ surface. The first is to discretely approximate metrics with random surfaces generated by gluing together regular polygons. Integrating a particular function uncovers of a generating function that obeys an infinite sequence of partial differential equations called the \emph{Korteweg-de Vries hierarchy}\index{KdV hierarchy}. The second method is to use supersymmetry to reduce this infinite dimensional integral to one over the finite dimensional submanifold of conformal metrics and hence to intersection numbers on the Deligne-Mumford compactification of $\mathcal{M}(S_{g,n})$. The (mathematically unjustified) belief that gravity is unique led to the conjecture that a certain generating function for these intersection numbers is a solution to the KdV hierarchy.

There is an established and varied family of proofs for Witten's conjecture, including:
\begin{itemize}
\item 
Kontsevich's proof via a fatgraph/ribbon graph model of a ``fattened" moduli space and matrix models techniques \cite{kontsevich};
\item
Okounkov-Pandharipande's proof via Hurwitz numbers, graph-based enumeration of branched covers of $\mathbb{C}P^1$ and matrix models \cite{op};
\item
Kazarian-Lando's proof by using the ELSV formula to relate Hurwitz numbers to intersection numbers on moduli space \cite{kl}; 
\item
Kim-Liu's proof using localization to obtain relations between Hodge integrals, which in turn imply Witten's conjecture \cite{kimliu}, and 
\item
Mirzakhani's proof via the Duistermaat-Heckman theorem and her Weil-Petersson volume recursion formula.
\end{itemize}

\subsubsection{Background.}
Our hitherto treatment of moduli space theory has been based on hyperbolic surfaces, however it is usual to consider $\mathcal{M}(S_{g,n})$ as the moduli space of genus $g$ Riemann surfaces with $n$ labeled marked points when defining $\psi$-classes\index{$\psi$-classes}.

A \emph{nodal Riemann surface}\index{nodal Riemann surface} is a connected complex space where every point is locally modeled upon either a disk in $\mathbb{C}$ or the neighborhood of a \emph{node}
\begin{align*}
\{(z,w)\in\mathbb{C}^2\mid zw=0\text{ and }|w|,|z|<1\}.
\end{align*}
The genus of a nodal Riemann surface is obtained by replacing each nodal neighborhood with an annulus and calculating the genus of the resulting topological surface. We call a nodal Riemann surface \emph{stable} if each connected component left after removing its nodes has negative Euler characteristic. Using uniformization, we may endow each of these components with a canonical hyperbolic metric with cusps at the nodes. Thus, stable nodal Riemann surfaces may be thought of as hyperbolic surfaces with length $0$ (cusp) interior ``geodesics".

\begin{figure}[h]
\begin{center}
\includegraphics[scale=0.3]{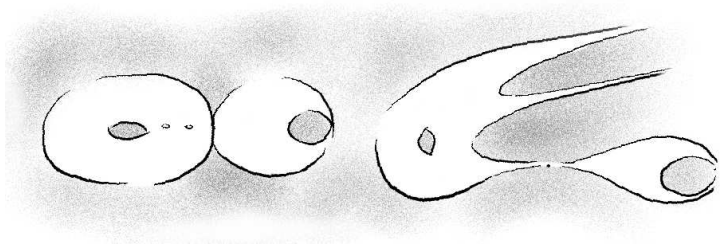}
\end{center}
\caption{A stable Riemann surface before and after uniformization.}
\label{pg21}
\end{figure}

The \emph{Deligne-Mumford compactification}\index{moduli space!Deligne-Mumford compactification} $\overline{\mathcal{M}}(S_{g,n})$ of $\mathcal{M}(S_{g,n})$ is the moduli space of genus $g$ stable nodal Riemann surfaces with $n$ labeled marked points. This is a smooth compactification of $\mathcal{M}(S_{g,n})$, and Wolpert \cite{wolpert1} showed that the Weil-Petersson form $\omega_{S_{g,n}}$ extends smoothly to $\overline{\mathcal{M}}(S_{g,n})$. From the hyperbolic perspective, $\overline{\mathcal{M}}(S_{g,n})$ is the moduli space of genus $g$ hyperbolic surfaces with $n$ labeled boundary cusps and possibly with some internal cusps. 

\begin{remark}
The hyperbolic viewpoint for $\overline{\mathcal{M}}(S_{g,n})$ admits a natural generalization in the way of moduli spaces $\overline{\mathcal{M}}(S_{g,n},\vec{b})$ of genus $g$ hyperbolic surfaces with $n$ labeled borders $\beta_1,\ldots,\beta_n$ of lengths $b_1,\ldots,b_n$, whereby we add in a compactification locus whose points correspond to bordered hyperbolic surfaces with interior cusps.
\end{remark}

A neighborhood $U$ around a non-orbifold point $[X]$ in $\overline{\mathcal{M}}(S_{g.n})$ may be thought of as a parameter space for surface structure deformations on $[X]$. Thus, we can build a surface bundle over the non-orbifold locus of $\overline{\mathcal{M}}(S_{g,n})$ so that the fiber over $[X]\in\overline{\mathcal{M}}(S_{g.n})$ \emph{is} the surface $X$ itself. In order to \emph{continuously} extend this construction over an orbifold point $[Y]$, the fiber over $[Y]$ must be $Y/\mathrm{Aut}(Y)$. To see this, consider a hyperelliptic surface $Y$ with involution $\iota:Y\rightarrow Y$. A path approaching a point $y$ on the fiber over $[Y]$ from nearby fibers is indistinguishable to a path approaching $\iota(y)$. Thus, for limits to exist, we must identify $y$ and $\iota(y)$. 

The resulting extended ``orbifold surface bundle":
\begin{align*}
P:\overline{\mathcal{C}}(S_{g,n})\rightarrow\overline{\mathcal{M}}(S_{g,n}),
\end{align*}
is referred to as the \emph{universal curve}\index{universal curve} over $\overline{\mathcal{M}}(S_{g,n})$.

Both the DM compactified moduli space $\overline{\mathcal{M}}(S_{g,n})$ and the space $\overline{\mathcal{C}}(S_{g,n})$ are naturally endowed with complex structures \cite{hm}, and there are $n$ holomorphic sections:
\begin{align*}
s_i:\overline{\mathcal{M}}(S_{g,n})&\rightarrow\overline{\mathcal{C}}(S_{g,n})\\
[X]&\to ([X],\beta_i),
\end{align*}
where $\beta_i$ denotes the $i$-th marked point on the fiber $X/\mathrm{Aut}(X)=P^{-1}([X])$ over $[X]$. Let $\mathcal{L}$ denote the holomorphic line bundle on $\overline{\mathcal{C}}(S_{g,n})$ that is contangent to the fibers of the universal curve. We define  $n$ \emph{tautological line bundles}\index{tautological line bundles} $\mathcal{L}_i:=s_i^*(\mathcal{L})$ on $\overline{\mathcal{M}}(S_{g,n})$ by pullback. The fiber of $\mathcal{L}_i$ over $[X]\in\overline{\mathcal{M}}(S_{g,n})$ is the cotangent space at the $i$-th marked point on $X$. Denote the first Chern class of $\mathcal{L}_i$ by $\psi_i:=\mathrm{c}_1(\mathcal{L}_i)$, they are called \emph{$\psi$-classes}\index{$\psi$-classes}.

\subsubsection{Statement of Witten's conjecture.}
Consider a collection of non-negative integers $a_1,\ldots,a_n\in\mathbb{Z}_{\geq0}$, define
\begin{align}
\langle \tau_{a_1},\ldots,\tau_{a_n}\rangle_g:=\left\{
\begin{array}{r l}
\int_{\overline{\mathcal{M}}(S_{g,n})}\prod_{i=1}^n\psi_i^{a_i},&\text{ if }\sum a_i=3g-3+n\\
0, &\text{ if }\sum a_i\neq 3g-3+n.
\end{array}
\right.
\end{align}
These are intersection numbers of $\psi$-classes, and the fact that $\psi_i$ is an even degree cohomology classes means that they commute and the ordering of the $a_i$ is unimportant. Note that given $a_1,\ldots,a_n$, there is at most one possible $g$ for which $\langle \tau_{a_1},\ldots,\tau_{a_n}\rangle_g$ might be non-zero (i.e.: when $\frac{1}{3}(\sum a_i+n)+1$ is a non-negative integer). We shall suppress writing the subscript $g$ when convenient. 

Consider the following generating function for these rational constants:
\begin{align}
F(t_0,t_1,\ldots):=\sum_{g=0}^\infty\; \sum_{a_1,\ldots,a_n\in\mathbb{Z}_{\geq0}}\langle\tau_{a_1},\ldots,\tau_{a_n}\rangle_g\prod_{r>0}\frac{t_r^{n_r}}{n_r!},
\end{align}
where the interior sum is taken over the set of finite multisets in $\mathbb{Z}_{\geq0}$ and $n_r$ is the number of times $r$ appears among $a_1,\ldots,a_n$.

\begin{conjecture}[Witten's conjecture]\index{Witten's conjecture}
The formal power series $\exp(F(\vec{t}))$ is annihilated by the infinite sequence $L_{-1},L_0,L_1,\ldots$ of differential operators given by:
\begin{align}
L_{-1}:=&-\frac{\partial}{\partial t_0}+\frac{1}{2}t_0^2+\sum_{i=0}^\infty t_{i+1}\frac{\partial}{\partial t_i},\\
L_0:=&-\frac{3}{2}\frac{\partial}{\partial t_1}+\frac{1}{16}+\frac{1}{2}\sum_{i=0}^\infty (2i+1)t_i\frac{\partial}{\partial t_i}\text{ and}\\
L_k:=&-\frac{(2k+3)!!}{2^{k+1}}\frac{\partial}{\partial t_k}+\sum_{i=0}^\infty \frac{(2k+2i+1)!!}{2^{k+1}\cdot (2i+1)!!} t_i\frac{\partial}{\partial t_{k+i}}\notag\\
&+\sum_{i+j=k-2}^k\frac{(2i+1)!!(2j+1)!!}{2^{k+2}}\frac{\partial^2}{\partial t_i\partial t_j},\text{ for }k\geq1.
\end{align}
\end{conjecture}

In terms of $\psi$-class intersection numbers, the equation $L_1(\exp F)=0$ being true is equivalent to the \emph{string equation}\index{string equation}
\begin{align}
\langle\tau_0,\tau_{a_1},\ldots,\tau_{a_n}\rangle_g=\sum_{a_i\neq0}\langle\tau_{a_1},\ldots,\tau_{a_i-1},\ldots,\tau_{a_n}\rangle
\end{align}
being true for any $g,n$ and $a_i$ such that $\sum a_i=3g-2+n$, in conjunction with the constraint that $\langle\tau_0,\tau_0,\tau_0\rangle_0=1$.

Similarly, the equation $L_0(\exp F)=0$ being true is equivalent to the \emph{dilaton equation}\index{dialon equation}
\begin{align}
\langle\tau_1,\tau_{a_1},\ldots,\tau_{a_n}\rangle_g=(2g-2+n)\langle\tau_{a_1},\ldots,\tau_{a_n}\rangle_g
\end{align}
being true for any $g,n$ and $a_i$ such that $\sum a_i=3g-3+n$, in conjunction with the constraint that $\langle\tau_1\rangle_1=\frac{1}{24}$.

For $a_1\geq1$, the equation $L_{a_1-1}(\exp F)=0$ is equivalent to the following condition on $\psi$-class intersection numbers
\begin{align}\label{kdvformula}
\langle\tau_{a_1},&\ldots,\tau_{a_n}\rangle_g\notag\\
&=
\sum_{j=2}^n\frac{(2a_1+2a_j-1)!!}{(2a_1+1)!!(2a_j-1)!!}\langle\tau_{a_2},\ldots,\tau_{a_1+a_j-1},\ldots,\tau_{a_n}\rangle_g\notag\\
&+\sum_{j+k=a_1-2}\frac{(2j+1)!!(2k+1)!!)}{2\cdot(2a_1+1)!!}\langle\tau_j,\tau_k,\tau_{a_2},\ldots,\tau_{a_n}\rangle_{(g-1)}\notag\\
&+\sum_{j+k=a_1-2}\frac{(2j+1)!!(2k+1)!!)}{2\cdot(2a_1+1)!!}\sum_{I_1\sqcup I_2\subset\{2,\ldots,n\}}
\langle\tau_j,\tau_{\vec{a}_{I_1}}\rangle \cdot\langle\tau_k,\tau_{\vec{a}_{I_2}}\rangle,
\end{align}
where $(2i+1)!!=1\times3\times\cdots\times(2i+1)$, the multiset $\vec{a}_I$ consists of every $a_i$ with $i\in I$ and we have suppressed the genus in the last line out of convenience.

\begin{remark}\label{onlyneed}
The operators $\{L_i\}$ generate a subalgebra of the Virasoro algebra\index{Virasoro algebra} with central charge $c=0$, and satisfy the following relation:
\begin{align}
[L_j,L_k]=(j-k)\cdot L_{j+k}.
\end{align}
Thus, for the purposes for proving Witten's conjecture, we only really need to verify the string equation, the dilaton equation and the $k=2$ case of \eqref{kdvformula}.
\end{remark}

\subsubsection{Idea of proof.}
We break down Mirzakhani's proof of Witten's conjecture into two steps:
\begin{enumerate}
\item
express the coefficients of the top degree monomials in $V_{g,n}(\vec{b})$ in terms of $\psi$-class intersection numbers over $\overline{\mathcal{M}}(S_{g,n})$,
\item
convert the coefficients of the top degree terms in Mirzakhani's volume recursion formula \eqref{recursionformula} into $\psi$-class intersection numbers, thus deducing the string equation, the dilaton equation and equation~\eqref{kdvformula} for all $k\geq1$.
\end{enumerate}
\vspace{1.5em}
\noindent\textbf{Step 1:}
The following theorem allows us to relate the volume of moduli space with $\psi$ (and $\kappa_1=\frac{[\omega]}{2\pi}$) class intersection numbers.
\begin{theorem}[{\cite[Thm.~4.4]{mirz2}}]\label{dictionarythm}
The coefficients of the volume polynomial
\begin{align}
V_{g,n}(\vec{b})=\sum_{\vec{a}:|\vec{a}|\leq3g-3+n} \int_{\overline{\mathcal{M}}(S_{g,n})}\frac{\psi_1^{a_1}\cdots\psi_n^{a_n}\omega^{3g-3+n-|\vec{a}|}}{2^{|\vec{a}|}\cdot\vec{a}!\cdot(3g-3+n-|\vec{a}|)!}
\cdot\vec{b}^{2\vec{a}},
\end{align}
where $|\vec{a}|:=a_1+\cdots+a_n$, $\vec{a}!:=a_1!\cdots a_n!$ and $\vec{b}^{2\vec{a}}:=b_1^{2a_1}\cdots b_n^{2a_n}$.
\end{theorem}

Much like how $\mathcal{L}_i$ is constructed, construct another line bundle $\mathcal{K}_i$ where the fiber over $[X]\in\overline{\mathcal{M}}(S_{g,n})$ is the tangent space of $X$ at its $i$-th marked point. The duality between $\mathcal{L}_i$ and $\mathcal{K}_i$ means that there is an orientation-reversing isomorphism between their associated (principal) circle bundles. Let $\mathcal{S}_i$ denote the circle bundle associated with $\mathcal{K}_i$, but with the opposite orientation. Then,
\begin{align*}
\mathrm{c}_1(\mathcal{S}_i)=-c_1(\mathcal{K}_i)=\mathrm{c}_1(\mathcal{L}_i)=\psi_i.
\end{align*}
Any fiber of $\mathcal{S}_i$ over a point $[X]\in\overline{\mathcal{M}}(S_{g,n})$ may be thought of as the set of (positive) directions going through the $i$-th marked point on $[X]$. When considered from the hyperbolic perspective, this is the same as the set of geodesic rays emanating from the $i$-th cusp and hence is in natural bijection with points on any embedded horocycle at the $i$-th cusp. This suggests another description of the total space $E(\mathcal{S}_i)$ of $\mathcal{S}_i$ as:
\begin{align*}
E(\mathcal{S}_i)=\left\{ ([X],p) \mid [X]\in\overline{\mathcal{M}}(S_{g,n})\text{ and }p\in\tilde{\beta}_i\right\},
\end{align*}
where $\tilde{\beta}_i$ is the length $\frac{1}{4}$ horocycle at the $i$-th cusp. 

We can also define such circle bundles $\mathcal{S}_i(\vec{b})$ over $\overline{\mathcal{M}}(S_{g,n},\vec{b})$ by setting $\tilde{\beta}_i$ to be the \emph{hypercycl}e consisting of points of distance $\sinh^{-1}(\frac{1}{\sinh(b_i/2)})$ from boundary $\beta_i$. The fiber product of all of the $\mathcal{S}_i$ is a $n$-torus bundle, and its total space is a leaf sitting in of the following moduli space:
\begin{align*}
\widehat{\mathcal{M}}(S_{g,n}):=\left\{([X],p_1,\ldots,p_n)\mid [X]\in\overline{\mathcal{M}}(S_{g,n},\vec{b}), \vec{b}\in\mathbb{R}_{>0}^n\text{ and } p_i\in\tilde{\beta}_i\right\}.
\end{align*}
By capping each boundary component of $S_{g,n}$ with a pair of pants of boundary lengths $(0,0,b_i)$, we obtain a genus $g$ hyperbolic surface $S_{g,2n}$ with $2n$ cusps and interior geodesics $\Gamma=(\gamma_1,\ldots,\gamma_n)$ that used to be $(\beta_1,\ldots,\beta_n)$ on $S_{g,n}$. It is easy to see that $\widehat{\mathcal{M}}(S_{g,n})=\overline{\mathcal{M}}(S_{g,n},\Gamma)$, and inherits a Weil-Petersson symplectic form $\omega_{S_{g,n},\Gamma}$ from this identification.

We see therefore that we have a Hamiltonian $T^n$-space $(\overline{\mathcal{M}}(S_{g,n},\Gamma),\omega_{S_{g,n},\Gamma})$, where the torus action twists the pairs of pants separated off by the $\beta_i$, and where the moment map is given by:
\begin{align*}
\tfrac{1}{2}\ell^2:=\left(\tfrac{1}{2}\ell_{\gamma_1}^2,\ldots,\tfrac{1}{2}\ell_{\gamma_n}^2\right).
\end{align*}
The Duistermaat-Heckman theorem \cite[Thm~2.5]{guillemin}, combined with the fact that the reduced space $\ell^{-1}(\vec{b})/T^n$ is symplectomorphic to $\overline{\mathcal{M}}(S_{g,n},\vec{b})$ implies Theorem~\ref{dictionarythm}.\\[1.5em]
\noindent\textbf{Step 2:} Let us denote the coefficient of $\vec{b}^{2\vec{a}}$ in $V_{g,n}(\vec{b})$ by $C_g(\vec{a})$. Observe that when $|\vec{a}|=3g-3+n$ (i.e.: when $|\vec{a}|$ is maximal), $C_g(\vec{a})$ is the coefficient of a top degree monomial in $V_{g,n}(\vec{b})$ and satisfies:
\begin{align}\label{tauidentity}
C_g(\vec{a})=\tfrac{1}{2^{|\vec{a}|}\cdot \vec{a}!}\langle \tau_{a_1},\tau_{a_2},\ldots,\tau_{a_n}\rangle_g.
\end{align}
Just as we may suppress the subscript $g$ when writing $\langle \tau_{a_1},\ldots,\tau_{a_n}\rangle$, we also suppress the $g$ in $C_g(\vec{a})$ when convenient.

Restrict Mirzakahni's volume recursion formula \eqref{recursionformula} to the coefficients of $\vec{b}^{2\vec{a}}$, where $|\vec{a}|=3g-3+n$, we obtain that:
\begin{align}\label{coeffrecursion}
(2a_1+1)&C_g(\vec{a})\notag\\
=&\sum_{j=1}^n\frac{(2a_1+2a_j-1)!}{(2a_1)!\cdot (2a_j)!}C_g(a_2,\ldots, a_1+a_j-1,\ldots,a_n)\\
+&\sum_{j+k=a_1-2}\frac{(2j+1)!(2k+1)!}{2\cdot (2a_1)!}C_{g-1}(j,k,a_2,\ldots,a_n)\notag\\
+&\sum_{j+k=a_1-2}\frac{(2j+1)!(2k+1)!}{2\cdot (2a_1)!}\sum_{I_1\sqcup I_2=\{2,\ldots,n\}} C(j,\vec{a}_{I_1})\cdot C(k,\vec{a}_{I_2}).\notag
\end{align}
Invoking \eqref{tauidentity} to translate these coefficients into intersection numbers, we immediately obtain \eqref{kdvformula} for $a_1\geq2$. Note that when $a_1=0,1$, we discard the last two lines in \eqref{coeffrecursion}. The resulting identity when $a_1=0$ is the string equation. And when $a_1=1$, the resulting identity is the dilaton equation. By remark~\ref{onlyneed}, this suffices to prove Witten's conjecture.

\subsection{Simple length spectrum growth rate}\index{simple length spectrum}
McShane and Rivin \cite{mr} first showed that the number of simple closed geodesics on a one-cusped hyperbolic torus of length less than $L$ has order $L^2$ growth. Rivin \cite{rivin} extended this result to general (orientable) hyperbolic surfaces with genus $g$ and $n$ boundary components to show that the growth rate is always of order $L^{6g-6+2n}$. In \cite{mirz3}, Mirzakhani refined these results by employing her volume computation to relate the precise rate of the asymptotic simple geodesic growth rate with the Thurston measure of certain open balls in the measured lamination space $\mathcal{ML}(S)$ of $S$. 

We introduce a little notation in order to state her main result. 

\subsubsection{Statement of the main result.} Let $\gamma$ be a simple closed geodesic on a bordered hyperbolic surface $S$ with genus $g$ and $n$ labeled geodesic borders of length $\vec{b}$. A simple geodesic may be regarded as an element  of the space $\mathcal{ML}(S)$ of \emph{measured (geodesic) laminations}\index{measured lamination} on $S$, and we can define a length function $\ell_\gamma:\mathcal{T}(S,\vec{b})\rightarrow\mathbb{R}_{>0}$. More generally, any measured geodesic lamination $\lambda$ is the limit of some sequence $\{a_i\gamma_i\}$ of weighted simple closed geodesics, and there is a well-defined generalized length function for measured laminations via a limiting process $a_i\ell_{\gamma_i}\to\ell_\lambda$. The function
\begin{align*}
\ell_{(\cdot)}(\cdot):\mathcal{ML}(S)\times\mathcal{T}(S,\vec{b})\rightarrow\mathbb{R}_{>0}
\end{align*}
is continuous \cite[Prop.~3]{bonahon}. Define the open ``unit ball''
\begin{align*}
B_{[X,f]}:=
\left\{
\lambda\in\mathcal{ML}(S)\mid \ell_\lambda ([X,f])\leq 1
\right\}.
\end{align*}
Although $B_{[X,f]}\subset \mathcal{ML}(S)$ depends upon the map $f$, its \emph{Thurston measure} $\mu_{\mathrm{Th}}(B_{[X,f]})$ depends only on the geometry of $X$. Thus, we may define a function $B:\mathcal{M}(S,\vec{b})\rightarrow\mathbb{R}_{\geq0}$ by taking
\begin{align*}
B([X]):=\mu_{\mathrm{Th}}(B_{[X,f]})
\end{align*}
for any labeling-preserving homeomorphism $f:S\rightarrow X$. Since $B$ is defined in terms of $\ell$ --- a continuous function, it must itself be continuous. For $\epsilon>0$ small enough so that simple closed geodesics of length shorter than $\epsilon$ cannot intersect, Mirzakhani obtained the following bounds for $B$ on the $\epsilon$-thin part of moduli space: for constants $C_1,C_2>0$, 
\begin{align*}
C_1\cdot\prod_{\alpha:\ell_\alpha([X,f])\leq\epsilon}\frac{1}{\ell_\alpha|\log(\ell_\alpha)|}\leq B([X])\leq C_2\cdot\prod_{\alpha:\ell_\alpha([X,f])\leq\epsilon}\frac{1}{\ell_\alpha},
\end{align*}
where the products are taken over all simple closed curves on $[X,f]$ shorter than $\epsilon$. The left inequality shows that $B$ is proper, and the right inequality shows that $B$ is dominated by an integrable function and hence is itself integrable. We define the number
\begin{align}
b(S)=\int_{\mathcal{M}(S,\vec{b})}B([X])\;\Omega_{S,\vec{b}}.
\end{align}
\begin{theorem}[{\cite[Thm.~1.1, 1.2]{mirz3}}]\label{simplegrowth}\index{simple length spectrum!growth rate}
Given a simple closed geodesic $\gamma$ on a hyperbolic surface $X$ with genus $g$ and $n$ labeled borders of length $\vec{b}$, define 
\begin{align*}
s_X(L,\gamma):=\mathrm{Card}\left\{\alpha\in\mathrm{Mod}(X)\cdot\gamma \mid \alpha \text{ has length }\leq L \right\}.
\end{align*}
Then, $n_\gamma:\mathcal{M}(S,\vec{b})\rightarrow\mathbb{R}_{>0}$ satisfies
\begin{align*}
n_\gamma([X]):=\lim_{L\to\infty}\frac{s_X(L,\gamma)}{L^{6g-6+2n}}
\end{align*}
is a continuous proper function. In particular, it is equal to
\begin{align}\label{growthformula}
\frac{c(\gamma)}{b(S)}\cdot B([X]),
\end{align}
where
\begin{align}
c(\gamma):=\lim_{L\to\infty}\frac{1}{L^{6g-6+2n}}\int_{\mathcal{M}(S,\vec{b})}s_X(L,\gamma)\; \Omega_{S,\vec{b}}
\end{align}
is a positive constant that depends only on the mapping class of $\gamma$ in $S$.
\end{theorem}
The main way in which Mirzakhani uses her volume integration techniques is as follows: she adapts the preimage-weighting argument used to unwrap the volume integral of moduli space to an integral of a function over an intermediate moduli space to see that:
\begin{align}\label{curvecountint}
\int_{\mathcal{M}(S,\vec{b})}s_X(L,\gamma)\; \Omega_{S,\vec{b}}=&\int_{\mathcal{M}(S,\gamma,\vec{b})}\chi_{\ell_\gamma^{-1}(0,L]}\;\Omega_{S,\gamma,\vec{b}}=\int_0^L\mathrm{Vol}(\ell_\gamma^{-1}(x))\;\mathrm{d}x.
\end{align}
Equation~\eqref{levelset} tells us the level set $\ell_\gamma^{-1}(x)$ is the product of $\mathbb{R}/x\mathbb{Z}$ and some $6g-8-2n$-dimensional moduli space. By Theorem~\ref{polynomial}, its volume $\mathrm{Vol}(\ell_\gamma^{-1}(x))$ is a product of $x$ and a degree $6g-8-2n$ polynomial $p(x)$ in $x$ (we treat the $\vec{b}$ as constants). Thus, the integral in \eqref{curvecountint} is equal to
\begin{align*}
P(L):=\int_0^L x\cdot p(x)\;\mathrm{d}x,
\end{align*}
which yields an order $6g-6+2n$ polynomial $P(L)$ in $L$ (and $\vec{b}$). This suffices to prove that $c(\gamma)=\lim_{L\to\infty}\frac{P(L)}{L^{6g-6+2n}}$ is positive and hence $n_\gamma$ is also strictly positive.

\begin{remark}
Since $c(\gamma)$ and $b(S)$ are constants, $n_\gamma$ is continuous and proper because $B$ is a continuous and proper function.
\end{remark}

\subsubsection{Idea of proof.}
We break down Mirzakhani's proof of Theorem~\ref{simplegrowth} into the following three steps:
\begin{enumerate}
\item
express $\frac{s_X(L,\gamma)}{L^{6g-6+2n}}$ as some discrete measure evaluated on $B_{[X,f]}\subset\mathcal{ML}(S)$,
\item
relate weak limits of these discrete measures to the Thurston measure,
\item
show that all such weak limits agree, by evaluating on $B_{[X,f]}$.
\end{enumerate}
\vspace{1em}
\noindent\textbf{Step 1:} Note that positive numbers homeomorphically act on $\mathcal{ML}(S)$ by multiplication. Given $L\in\mathbb{R}_{>0}$ and a Borel subset $V\subset \mathcal{ML}(S)$, we define the following discrete measure:
\begin{align}
\mu_{L,\gamma}(V):=\frac{\mathrm{Card}\left(L\cdot V\cap \mathrm{Mod}(S)\cdot \gamma \right)}{L^{\mathrm{dim}_{\mathbb{R}}\mathcal{ML}(S)=6g-6+2n}},
\end{align}
and observe that
\begin{align}
\frac{s_X(L,\gamma)}{L^{6g-6+2n}}=\mu_{L,\gamma}(B_{[X,f]}).
\end{align}
Therefore, one way to study the asymptotic behavior of $s_X(L,\gamma)$ is to study the limiting behavior of the family $\{\mu_{L,\gamma}\}_{L\in\mathbb{R}}$ as $L$ tends to infinity.\\[1.5em]
\noindent\textbf{Step 2:} Let $\mathcal{ML}(S,\mathbb{Z})\subset\mathcal{ML}(S)$ denote the set of \emph{integral multi-curves} on $S$. Since $\mathcal{ML}(S,\mathbb{Z})$ is a lattice in $\mathcal{ML}(S)$, the family of measures defined by
\begin{align*}
\nu_{L,\gamma}(V):=\frac{\mathrm{Card}\left(L\cdot V\cap \mathcal{ML}(S,\mathbb{Z})\right)}{L^{6g-6+2n}}
\end{align*}
approximates the Thurston measure $\mu_{\mathrm{Th}}$ as $L$ tends to infinity. The fact that $\mu_{L,\gamma}(V)\leq\nu_{L,\gamma}(V)$ means that $\{\mu_{L,\gamma}\}$ is a bounded family and that any weak limit $\mu_{\vec{L},\gamma}$, with respect to some sequence of increasing indices $\vec{L}=(L_1,L_2,\ldots)$, must be absolutely continuous with respect to $\mu_{\mathrm{Th}}$. 

Since $\mu_{\vec{L},\gamma}$ and $\mu_{\mathrm{Th}}$ are both mapping class group invariant, and 
$\mathrm{Mod}(S)$ acts ergodically on $\mathcal{ML}(S)$, the Radon-Nikodym derivative $\frac{\mathrm{d}\mu}{\mathrm{d}\nu}$ must be equal to some constant $C_\vec{L}\in\mathbb{R}_{>0}$. Thus:
\begin{align*}
\mu_{\vec{L},\gamma}=C_\vec{L}\cdot\mu_{Th}.
\end{align*}
We finish showing that $C_\vec{L}$ is independent of $\vec{L}$.\\[1.5em]
\noindent\textbf{Step 3:} Just as $B([X]):=\mu_{\mathrm{Th}}(B_{[X,f]})$ defines a real function on $\mathcal{M}(S,\vec{b})$, we can define the following functions
\begin{align*}
B_{L,\gamma}([X]):=\mu_{L,\gamma}(B_{[X,f]})\text{ and }B_{\vec{L},\gamma}([X]):=\mu_{\vec{L},\gamma}(B_{[X,f]}).
\end{align*}
Integrating $B_{\vec{L},\gamma}$ over the moduli space $\mathcal{M}(S,\vec{b})$, we obtain:
\begin{align*}
\int_{\mathcal{M}(S,\vec{b})} B_{\vec{L},\gamma}([X])\; \Omega_{S,\vec{b}}
&=\int_{\mathcal{M}(S,\vec{b})}\lim_{i\to\infty}B_{L_i,\gamma}([X])\;\Omega_{S,\vec{b}}\\
\Rightarrow 
C_\vec{L}\cdot \int_{\mathcal{M}(S,\vec{b})} B([X])\; \Omega_{S,\vec{b}}
&=\lim_{L\to\infty}\int_{\mathcal{M}(S,\vec{b})}\frac{s_X(L,\gamma)}{L^{6g-6+2n}}\;\Omega_{S,\vec{b}}\\
\Rightarrow C_\vec{L}
&=\frac{c(\gamma)}{b(S)}.
\end{align*}
We see therefore that $C_\vec{L}$ is $\vec{L}$-independent, and
\begin{align}\label{measures}
\mu_{\vec{L},\gamma}=\frac{c(\gamma)}{b(S)}\cdot\mu_{\mathrm{Th}}.
\end{align}
Evaluating \eqref{measures} on $B_{[X,f]}$ yields equation~\eqref{growthformula}.

\subsubsection{Additional remarks.} We first note that the collection of simple closed geodesics on $S$ decomposes as a finite disjoint union of mapping class group orbits 
\begin{align*}
\mathrm{Mod}(S)\cdot\gamma_1\sqcup\ldots\sqcup\mathrm{Mod}(S)\cdot\gamma_k,
\end{align*}
and the asymptotic growth rate of the number $s_X(L)$ of simple closed geodesics on $X$ of length less than $L$ is given by:
\begin{align*}
s_X(L)\sim \frac{B([X])}{b(S)}\cdot\left(\sum_{i=1}^k c(\gamma_i)\right) L^{\mathrm{dim}_{\mathbb{R}}\mathcal{ML}(S)}.
\end{align*}
Moreover, Mirzakhani's proof applies to rational multicurves in general and not just to simple closed geodesics.

Mirzakhani asserts that equation~\eqref{growthformula} holds for closed curves in general. One possible way to prove this might be to use trace relations or hyperbolic trigonometry to show that when the length of a closed curve $\gamma\subset S$ on a marked hyperbolic surface $[X,f]$ is long, it is roughly ``linear'' in the lengths of certain disjoint simple closed geodesics $\{\gamma_i\}$. Then, the asymptotic growth rate for $\gamma$ is given by the growth rate of a specific multi-curve built from these $\{\gamma_i\}$. In fact, this idea applies rather easily to let us obtain asymptotic growth rates of orthogeodesics (and related geometric objects). 

More importantly, Mirzakhani also claims that \eqref{growthformula} holds for \emph{geodesic currents} --- a far-reachingly general class of limiting objects built from closed geodesics in much the same way that measured geodesic laminations are built from simple closed curves --- with potential applications in studying the length spectrum growth rates for non-positively curved surfaces and Hitchin representations or perhaps the growth rate of geometric intersection numbers of pairs of closed curves\ldots etc.

As a final remark, it is curious to note that polynomial growth rates do not (always) hold in the non-orientable context. In particular, Norbury and the author have shown that the order of the growth rate for 1-sided simple closed geodesics on any thrice-cusped projective plane is between $O(L^{2.430})$ and $O(L^{2.477})$ \cite{hn}\index{simple length spectrum!non-orientable surface growth rate}.

\appendix
\section{$F_{2k+1}(b)$}\label{appendixsec}
We demonstrate how to explicitly compute the integral $F_{2k+1}(b)$ defined in Subsection~\ref{polynomialsubsec}. 
\begin{align*}
F_{2k+1}(b)=
&\int_0^\infty \left(\frac{x^{2k+1}}{1+\exp(x+b)}+\frac{x^{2k+1}}{1+\exp(x-b)}\right)\mathrm{d}x\\
=&\int_b^\infty\frac{(x-b)^{2k+1}\;\mathrm{d}x}{1+\exp(x)} + \int_{-b}^\infty\frac{(x+b)^{2k+1}\;\mathrm{d}x}{1+\exp(x)}\\
=&\int_{-b}^0\frac{(x+b)^{2k+1}\;\mathrm{d}x}{1+\exp(x)} -\int_0^b\frac{(x-b)^{2k+1}\;\mathrm{d}x}{1+\exp(x)}\\
&+\int_0^\infty\left(\frac{(x+b)^{2k+1}+(x-b)^{2k+1}}{1+\exp(x)}\right)\mathrm{d}x\\
=&\int_{-b}^0(x+b)^{2k+1}\left(\frac{1}{1+\exp(x)}+\frac{1}{1+\exp(-x)}\right) \;\mathrm{d}x\\
&+\int_0^\infty\left(\frac{2\sum_{i=1}^{k+1} {2k+1 \choose 2i-1} \cdot b^{2k+2-2i}\cdot x^{2i-1}}{1+\exp(x)}\right)\mathrm{d}x\\
=&\frac{b^{2k+2}}{2k+2}+2\sum_{i=1}^{k+1}b^{2k+2-2i}{2k+1\choose 2i-1}\int_0^\infty\frac{x^{2i-1}\;\mathrm{d}x}{1+\exp(x)}.
\end{align*}
To compute $\int_0^\infty \frac{x^{2i-1}\mathrm{d}x}{1+\exp(x)}$, we use the fact that
\begin{align*}
\frac{1}{1+\exp(x)}=\frac{1}{\exp(x)-1}-\frac{2}{\exp(2x)-1}
\end{align*}
to see that
\begin{align*}
\int_0^\infty \frac{x^{2i-1}\mathrm{d}x}{1+\exp(x)}=&\int_0^\infty \frac{x^{2i-1}\mathrm{d}x}{\exp(x)-1}+\int_0^\infty \frac{2x^{2i-1}\mathrm{d}x}{\exp(2x)-1}\\
=&\int_0^\infty \frac{x^{2i-1}\mathrm{d}x}{\exp(x)-1}+\int_0^\infty \frac{2^{1-2i} y^{2i-1}\mathrm{d}y}{\exp(y)-1}\\
=&\zeta(2i)(2i-1)!(1-2^{1-2i}).
\end{align*}
Since $\zeta(0)=-\frac{1}{2}$, we may express $F_{2k+1}(b)$ as the following polynomial
\begin{align*}
F_{2k+1}(b)={(2k+1)!}\sum_{i=0}^{k+1}\zeta(2i)(2^{2i+1}-4)\frac{b^{2k+2-2i}}{(2k+2-2i)!}.
\end{align*}


\frenchspacing
\begin{thebibliography}{1}

\bibitem{birmanseries}
J. S. Birman, C. Series, Geodesics with bounded intersection number on surfaces are sparsely distributed. {\it Topology}, 24 (1985), no. 2, 217--225.

\bibitem{bonahon}
F. Bonahon, The geometry of Teichm\"uller space via geodesic currents. {\it Invent.\ Math.}, 92 (1988), no. 1, 139--162.

\bibitem{buser}
P. Buser, {\it Geometry and spectra of compact Riemann surfaces}. Progress in Mathematics, 106, Birkh\"auser Boston Inc., Boston, MA, 1992.

\bibitem{guillemin}
V. Guillemin, {\it Moment maps and combinatorial invariants of Hamiltonian $T^n$-spaces}. Progress in Mathematics, 122, Birkh\"auser Boston, Inc., Boston, MA, 1994.

\bibitem{hm}
J. Harris and I. Morrison, {\it Moduli of curves}, Graduate Texts in Mathematics, 187, Springer-Verlag, New York, 1998.

\bibitem{hn}
Y. Huang and P. Norbury, Simple geodesics and Markoff quads. {\tt math.GT/} {\tt 1312.7089}.

\bibitem{kl}
M. E. Kazarian and S. K. Lando, An algebro-geometric proof of Witten's conjecture. {\it J. Amer. Math. Soc.}, 20 (2007), no. 4, 1079--1089.

\bibitem{kimliu}
Y.-S. Kim and K. Liu, A new approach to deriving recursion relations for the Gromov-Witten theory. In {\it Topology and physics}, Nankai Tracts Math., 12, World Sci.\ Publ., Hackensack, NJ, 2008, 195--219.

\bibitem{kontsevich}
M. Kontsevich, Intersection theory on the moduli space of curves and the matrix Airy function. {\it Comm.\ Math.\ Phys.}, 147 (1992), no. 1, 1--23.

\bibitem{primer}
B. Farb and D. Margalit, {\it A primer on mapping class groups}. Princeton Mathematical Series, 49, Princeton University Press, Princeton, NJ, 2012.

\bibitem{mcshane}
G. McShane, {\it A remarkable identity for lengths of curves}. Ph.D. Thesis, University of Warwick, 1991.

\bibitem{mr}
G. McShane and I. Rivin, A norm on homology of surfaces and counting simple geodesics. {\it Internat.\ Math.\ Res.\ Notices}, 1995, no. 2, 61--69.

\bibitem{mirz3}
M. Mirzakhani, Growth of the number of simple closed geodesics on hyperbolic surfaces. {\it Ann. of Math. (2)}, 168 (2008), no. 1, 97--125.

\bibitem{mirz1}
M. Mirzakhani, Simple geodesics and Weil-Petersson volumes of moduli spaces of bordered Riemann surfaces. {\it Invent. Math.}, 167 (2007), no. 1, 179--222.

\bibitem{mirz2}
M. Mirzakhani, Weil-Petersson volumes and intersection theory on the moduli space of curves. {\it J. Amer. Math. Soc.}, 20 (2007), no. 1, 1--23.

\bibitem{nn1}
M. N\"a\"at\"anen and T. Nakanishi, Weil-Petersson areas of the moduli spaces of tori. {\it Results Math.}, 33 (1998), no. 1-2, 120--133.

\bibitem{nn2}
T. Nakanishi and M. N\"a\"at\"anen, Areas of two-dimensional moduli spaces. {\it Proc. Amer. Math. Soc.}, 129 (2001), no. 11, 3241--3252.

\bibitem{op}
A. Okounkov and R. Pandharipande, Gromov-Witten theory, Hurwitz numbers, and matrix models. In {\it Algebraic geometry, Seattle 2005}, ({S}eattle, WA, 2005), Proc.\ Sympos.\ Pure Math., 80, Amer. Math. Soc., Providence, RI, 2009, 325--414.

\bibitem{parlier}
H. Parlier, Simple closed geodesics and the study of Teichm\"uller spaces. {\it Handbook of Teichm\"uller theory. Vol. IV}, IRMA Lect.\ Math.\ Theor.\ Phys., 19, Eur.\ Math.\ Soc., Z\"urich, 2014, 113--134.

\bibitem{penner}
 R. C. Penner, Weil-Petersson volumes. {\it J.\ Differential Geom.}, 35 (1992), no. 3, 559--608.

\bibitem{rivin}
I. Rivin, Simple curves on surfaces. {\it Geom. Dedicata}, 87 (2001), no. 1-3, 345--360.

\bibitem{teich}
O. Teichm\"uller, {\it Extremale quasikonforme Abbildungen und quadratische Differentiale}. Abh.\ Preuss.\ Akad.\ Wiss.\ Math.\ -Nat.\ Kl., 1939 (1940).

\bibitem{witten}
E. Witten, Two-dimensional gravity and intersection theory on moduli space. In {\it Surveys in differential geometry ({C}ambridge, {MA}, 1990)}, Lehigh Univ., Bethlehem, PA, 1991, 243--310.

\bibitem{wolpert1}
S. Wolpert, On the homology of the moduli space of stable curves. {\it Ann. of Math. (2)}, 118 (1983), no. 3, 491--523.

\bibitem{wolpert2}
S. Wolpert, On the K\"ahler form of the moduli space of once punctured tori. {\it Comment.\ Math.\ Helv.}, 58 (1983), no. 2, 246--256.

\bibitem{wolpert3}
S. Wolpert, On the Weil-Petersson geometry of the moduli space of curves. {\it Amer.\ J.\ Math.}, 107 (1985), no. 4, 969--997.

\bibitem{zograf2}
P. G. Zograf, The Weil-Petersson volume of moduli spaces of curves of small genus. {\it Funktsional.\ Anal.\ i Prilozhen.}, 32 (1998), no. 4, 78--81.

\bibitem{zograf1}
P. G. Zograf, The Weil-Petersson volume of the moduli space of punctured spheres. In {\it Mapping class groups and moduli spaces of {R}iemann surfaces} (G\"ottingen, 1991/Seattle, WA, 1991), Contemp.\ Math., 150, Amer.\ Math.\ Soc., Providence, RI, 1993, 367--372.


\end{thebibliography}
\end{document}